\documentclass[12pt]{article}
\usepackage{amsmath,amssymb,epsfig}

\newcommand{\be}{\begin{equation}}
\newcommand{\ee}{\end{equation}}
\def\H{{\cal H}}
\def\B{{\cal B}}
\def\N{{\cal N}}
\newcounter{unnumber}

\newtheorem{thrm}{Theorem}

\newtheorem{prop}{Proposition}

\newtheorem{de}{Definition}
\newtheorem{prf}[unnumber]{Proof}
\newenvironment{pf}{\begin{prf} \rm}{\end{prf}}

\numberwithin{equation}{section}

\begin{document}

\thispagestyle{empty} \phantom{}\vspace{-1.5cm}
\rightline{\today} \vspace{-0.5cm}

\begin{center}
  {\bf \LARGE Representations of the
ultrahyperbolic BMS group $\H \B
$.
\newline
II.  Determination of the representations induced from
infinite
little groups
}
\end{center}




\vspace{1.3cm} \centerline{
Evangelos Melas $^{a,}$ \footnote{emelas@econ.uoa.gr $\&$
evangelosmelas@yahoo.co.uk }} \vspace{.6cm}
{\em \centerline{$^a$ Department of Economics, Unit of Mathematics and Informatics,}
\centerline{ Athens, Greece}}



\vspace{1.0truecm}
\begin{abstract}
The ordinary Bondi$-$Metzner$-$Sachs (BMS) group $B$ is the common
asymptotic symmetry group of all asymptotically flat Lorentzian
space-times. \ As such, $B$ is the best candidate for the
universal symmetry group of General Relativity. \ However, in
studying quantum gravity, space$-$times with signatures other than
the usual Lorentzian one, and complex space$-$times, are frequently
considered. \ Generalisations of $B$ appropriate to these other
signatures have been defined earlier. In particular,
the generalisation
$\H \B$, a BMS group appropriate to
the ultrahyperbolic signature (+,+,$-$,$-$),
has been defined in a previous paper where it was shown that all the
strongly continuous unitary irreducible representations (IRs) of $\H \B $
can be obtained with the Wigner$-$Mackey's inducing method and that all the little groups of
$\H \B $ are compact.
\ Here we describe in detail all the infinite little groups of $\H \B$  and we find
the IRs of $\H \B$  induced by them.
\end{abstract}

\bigskip \bigskip
\newpage
\pagenumbering{arabic}









\section{Introduction
}
\indent

The best candidate for the universal symmetry group of General
Relativity (G.R), in any signature, is the so called
Bondi$-$Metzner$-$Sachs (BMS) group $B$. These groups have recently been
described \cite{mac1} for all possible signatures and all possible
complex versions of GR as well.

\noindent In earlier papers \cite{mac1,macmel,mel,mel1,mel2}  it has
been argued that the IRs of the BMS group and of its
generalizations in complex space$-$times as well as in space$-$times
with Euclidean or Ultrahyperbolic signature are what really lie
behind the full description of (unconstrained) moduli spaces of
gravitational instantons.
Kronheimer \cite{Kronheimer, Kronheimer1} has given a description
of these instanton moduli spaces for {\it Euclidean} instantons.
However, his description only partially describes the moduli
spaces, since it still involves {\it constraints}. Kronheimer does
not solve the constraint equations, but it has been argued
\cite{mac1,mel2} that IRs of BMS group (in the relevant signature)
give an {\it unconstrained} description of these same moduli
spaces.

The original BMS group $B$ was discovered by Bondi,
Metzner and Van der Burg \cite{Bondi}
for asymptotically flat space$-$times which were axisymmetric, and
by Sachs \cite{Sachs1} for general asymptotically flat
space$-$times, in the usual Lorentzian signature. The group $\mathcal H \mathcal B$
is a different generalised BMS group, namely  one appropriate
to the `ultrahyperbolic' signature, and asymptotic flatness in
null directions introduced in \cite{Mel}.

Recall that the ultrahyperbolic version of
Minkowski space is the vector space $R^{4}$ of row vectors with 4
real components, with scalar product defined as follows. \ Let
$x,y\in R^{4}$ have components $x^{\mu }$ and $y^{\mu }$
respectively, where $\mu =0,1,2,3$. \ Define the scalar product
$x.y$ between $x$ and $y$ by
\begin{equation}
x.y=x^{0}y^{0}+x^{2}y^{2}-x^{1}y^{1}-x^{3}y^{3}.
\end{equation}
Then the ultrahyperbolic version of Minkowski space, sometimes
written $R^{2,2}$, is just $R^{4}$ with this scalar product.

In \cite{Mel} it was shown that
\begin{thrm}
\label{omorpho} The group $\H\B
$   can be realised as \be
\H\B
= L^{2}(\mathcal P,\lambda ,R)
\bigcirc\!\!\;\!\!\!\;\!\!\!\!s \  _{T} G^{2} \ee with semi$-$direct
product specified by \be
\label{an}
(T(g,h)\alpha)(x,y)=k_{g}(x)s_{g}(x)k_{h}(w)s_{h}(w)\alpha
(xg,yh), \ee
\end{thrm}
where $\alpha \in L^{2}(\mathcal P,\lambda ,R)$ and $(x,y)\in \mathcal P$.
For ease of notation, we   write $\mathcal P$
for the torus ${\rm T} \simeq P_{1}(R)\times P_{1}(R)$, $P_{1}(R)$ is the one$-$dimensional real projective space,
and $\mathcal G$ for $G\times G$, $G=SL(2,R)$. In analogy to $B$, it is natural to choose a measure $\lambda$
on $\mathcal P$ which is invariant under
the maximal compact subgroup ${ \rm S}{\rm O}(2)\times {\rm S}{\rm
O}(2)$ of $\mathcal G $. $L^{2}(\mathcal P,\lambda ,R)$ is the
separable Hilbert space of real$-$valued functions defined on $\mathcal P$.




Moreover, if $g \in G$ is
\be
\left[ \begin{array}{lr}
   a    & b \\
   c    & d
\end{array}  \right],
\ee then the components $ x_{1},x_{2} $ of ${\bf x} \in R^{2}$ transform
linearly, so that the ratio $ x= {x_{1}}/{x_{2}} $ \  transforms
fraction linearly. Writing $ xg $ \ for the transformed ratio, \be
\label{aliki} xg=
\frac {({\bf x}g)_{1}}{({\bf x}g)_{2}}=\frac
{x_{1}a + x_{2}c} {x_{1}b+x_{2}d}=
\frac {xa+c}{xb+d}. \ee
 The factors $k_{g}(x)$ and $s_{g}(x)$ on the right hand side of (\ref{an}) are defined by
\be
\label{donkey}
k_{g}(x)=
\left\{ \frac {(xb+d)^{2}+(xa+c)^{2}}{1+x^{2}} \right\}
^{\frac {1}{2}},
\ee

\be
\label{mule}
s_{g}(x)=
\frac {xb+d}{|xb+d|},
\ee
 with similar formulae for $yh$, $k_{h}(y)$ \ and $s_{h}(y)$.

\indent
It is well known that the topological dual of a Hilbert space can be identified
with the Hilbert space itself, so that we
have $ {L^{2}}^{'}(\mathcal P,\lambda ,R) \simeq
L^{2}(\mathcal P,\lambda ,R).$
~ In fact, given a continuous linear functional
 $\phi \in {L^{2}}^{'}(\mathcal P,\lambda ,R)$, we can write,
 for $\alpha \in L^{2}(\mathcal P,\lambda ,R)$
\be
(\phi,\alpha)=<\phi,\alpha>
\ee
where the function $\phi \in L^{2}(\mathcal P,\lambda ,R)$ on the right is uniquely determined
by (and denoted by the same symbol as) the linear functional $\phi \in {L^{2}}^{'}(\mathcal P,\lambda ,R)
$ on the left.
The representation theory of
$\H\B
$
is governed by the dual
action $T'$ of $\mathcal G$ on the topological dual
$ {L^{2}}^{'}(\mathcal P,\lambda ,R)$ of $L^{2}(\mathcal P,\lambda ,R).$
The dual action $T'$ is defined by:
\be <T'(g,h)\phi,\alpha>=<\phi,T(g^{-1},h^{-1})\alpha>\cdot \ee A
short calculation
gives
\be \label{tzatziki}
(T'(g,h)\phi)(x,y)=k_{g}^{-3}(x)s_{g}(x)k_{h}^{-3}(y)s_{h}(y)\phi(xg,yh).
\ee
Now, this action $T'$ of $\mathcal G$
on ${L^{2}}^{'}(\mathcal P,\lambda ,R)$,
given explicitly above, is like the action $T$ of
$\mathcal G$ on $L^{2}(\mathcal P,\lambda ,R)$, continuous. The `little group' $L_{\phi}$
of any  $\phi \in {L^{2}}^{'}(\mathcal P,\lambda ,R)$ is the stabilizer \be
L_{\phi}=\{(g,h)\in \mathcal G
\ | \ T'(g,h)\phi=\phi\}. \ee By
continuity, $L_{\phi}\subset \mathcal G$ is a closed subgroup.

Attention is confined to measures on
${L^{2}}^{'}(\mathcal P,\lambda ,R)$ which are concentrated on
single orbits of the $
\mathcal G
-$action $T^{\prime }.$ These measures
give rise to IRs of $ \H \B $ which are induced in a sense
generalising \cite{mac6} Mackey's \cite{Wigner,Mackey,Mackey1,Simms,Isham,Mackey2}.
In fact in \cite{Mel} it was shown that \it all \normalfont
IRs of the $\H \B$
with the Hilbert topology are
derivable by the inducing construction.
The inducing construction
is realized as follows. Let $ \mathcal O \subset
{L^{2}}^{'}(\mathcal P,\lambda ,R)$ be any orbit of the dual action
${\;T}^{\prime }$ of
$
\mathcal G
$ on ${L^{2}}^{'}(\mathcal P,\lambda ,R).$
There is a natural homomorphism $ \mathcal O = \mathcal G \phi_{\rm{o}} \simeq
\mathcal G
/L_{\phi_{\rm{o}}} $ where $L_{\phi_{\rm{o}}} $ is the `little group' of the point
$\phi_{\rm{o}} \subset \mathcal O .$ Let $U$ be a continuous
irreducible unitary representation of $L_{\phi_{\rm{o}}}$ on a Hilbert space
$D.$
Every coset space  $\mathcal G / L_{\phi_{\rm{o}}}
$ can be  equipped with a unique
class of measures which are quasi$-$invariant under the  action
${\;T}$ of $G^{2}$. Let $ \mu $ be any one  of these. Let $
D_{\mu}= L^{2}(\mathcal G/ L_{\phi_{\rm{o}}} , \mu, D)$ be the
Hilbert space of functions $  f :  \mathcal G/ L_{\phi_{\rm{o}}}  \rightarrow D$
which are square integrable with respect to $\mu.$ From a
given $\phi_{\rm{o}} $ and any continuous irreducible unitary
representation $U$ of $L_{\phi_{\rm{o}}}$ on a Hilbert space $D$
a continuous irreducible unitary representation of $ \H \B \ $
on $D_{\mu}$ can be constructed. The representation is said to
be induced from $U$ and $\phi_{\rm{o}}$.
Different points of an orbit $\mathcal G \phi$
have conjugate little groups and give rise to equivalent
representations of $ \H \B $.

To conclude,
every irreducible representation of
$\H \B$ is obtained by the inducing construction for each $\phi_{\rm{o}} \in  {L^{2}}^{'}(\mathcal P,\lambda ,R)$
and each irreducible representation $U$ of $L_{\phi_{\rm{o}}}$.
All the little groups $L_{\phi_{\rm{o}}}$ of $ \H \B $
are compact and they up to conjugation subgroups of
$ \rm {S O(2)} \times  \rm { S O (2)} $. They
include groups
which are finite as well as groups which are
infinite, both  connected
and not$-$connected \cite{macmel}.
Therefore the construction of the IRs of $\mathcal H \mathcal B$
involves at the first instance the classification of all the subgroups of
$ \rm {S O(2)} \times  \rm { S O (2)} $.

In particular the problem of
constructing the IRs of $\mathcal H \mathcal B$
induced from {\it finite} little groups
reduces to a seemingly very simple task; that of classifying all
subgroups of the Cartesian product group $\; {\rm C}_{{\rm n}}
\times {\rm C}_{{\rm m}}, \;$ where $\;{\rm C}_{{\rm r}}\;$ is the
cyclic group of order $\; {\rm r}, \;$ $\; {\rm r}\;$ being
finite. Surprisingly, this task is less simple than it may appear
at first sight. It turns out \cite{mel} that the solution is constructed from
the `fundamental cases' $\; {\rm n}={\rm p}^{a}, \;$ $\; {\rm
m}={\rm p}^{\beta}, \;$ (n,m are powers of the same prime), via
the prime decomposition of m and n. The classification of all the
subgroups of $\; {\rm C}_{{\rm n}} \times {\rm C}_{{\rm m}}, \;$
was given in \cite{mel}. The explicit construction of the IRs of
$\mathcal H \mathcal B$ induced from finite little groups is not
a trivial task and it will be given in a different paper \cite{Mel1}.


The infinite not$-$connected
little groups
were given  in \cite{macmel}.
In this paper we concentrate on describing in detail all the
infinite connected little groups
and on finding the Hilbert spaces of the invariant functions
$\phi_{\rm{o}}$ of all the infinite little groups $L_{\phi_{\rm{o}}}$,
both connected and non$-$connected.
Finally we give explicitly the operators of the IRs  of
$\H \B$ induced from all the infinite little groups.



\indent
The paper is organised as follows: In Section \ref{s1} 
we construct the infinite connected subgroups of $ \rm  { S O(2)} \times \rm {S O (2)}. $ 
In Section \ref{s2} we find all the infinite potential little 
groups, both connected and non$-$connected. In Section \ref{s3} we prove 
that \it all \normalfont the infinite potential little groups are actual.
In Section \ref{s4} we construct explicitly the IRs of $\H \B$ induced 
from the infinite little groups. 
In Section \ref{s5} we make some remarks about the  IRs of $\H \B$ 
constructed in Section \ref{s4} by the inducing method.
 In Appendix \ref{s6} we give the elementary regions associated 
 with the action of the little groups on the torus $P_{1}(R) \times P_{1}(R).$ Finally 
In Appendix \ref{s7} we  construct 
$\;\mathcal G-{\rm quasi}-{\rm invariant}$ measures on the orbits $ 
\mathcal G \phi_{\rm{o}}$ of the action of the little groups $L_{\phi_{\rm{o}}} $ 
on the topological dual  ${L^{2}}^{'}(\mathcal P,\lambda ,R)$, measures which are necessary 
in order 
to give the operators of the induced IRs of $\H \B$ in Section \ref{s4}.

\section{Infinite Connected Subgroups of $ \mathbf  { S O(2)\times  S O (2)} $ }

\label{s1}

\indent

\ In \cite{Mel} it was shown that the little groups of $\H \B $ are compact subgroups of $ {\rm
S}{\rm L}(2,R) \times {\rm
S}{\rm L}(2,R)$.
The maximal compact subgroup of ${\rm
S}{\rm L}(2,R) \times {\rm
S}{\rm L}(2,R)$
is just the subgroup ${\rm S}{\rm O}(2)\times
{\rm S}{\rm O}(2)$. \ That is to say, if $H$ is a compact subgroup
of $G\times G$, then some conjugate $gHg^{-1}$ of $H$ is a compact
subgroup of ${\rm S}{\rm O}(2)\times {\rm S}{\rm O}(2)$. \ In the
representation theory of $ \; \H \B $,
the little groups
are only significant up to conjugacy. \ So, after a possible
conjugation, we may take every little group to be a compact (or
equivalently closed) subgroup of $ \; K={\rm S}{\rm O}(2)\times
{\rm S}{\rm O}(2). \;$

\begin{de}
For each closed subgroup $S$ of $ \; \mathcal G$,\ define its
maximal invariant (closed) vector space $ \; A(S)\subset {L^{2}}^{'}(\mathcal P,\lambda ,R) \;
$ by \be A(S)=\{\phi \in {L^{2}}^{'}(\mathcal P,\lambda ,R) \; | \; s\phi=\phi \quad for \ all
\quad s\in S \}. \ee
\end{de}
The little groups of \ $ \; \H\B
\; $\ may now be found by
following the programme proposed by McCarthy \cite{mac6}.
\vspace*{3mm}

{\em Programme}

(I) Determine, up to conjugacy, all compact subgroups $S \subset
{\rm S}{\rm O}(2)\times {\rm S}{\rm O}(2)$ .

(II) For each such $S$ find the vector space $A(S)$.

(III) Find the subgroups $S$ for which $A(S)$ contains elements
fixed under no properly larger group $S' \supset S$. These last
subgroups, and only these are the little groups.

\vspace{0.5cm} \noindent Now we proceed to carry out step (I) of
the programme.
The compact subgroups of $\; {\rm S}{\rm
O}(2)\times {\rm S}{\rm O}(2) \;$ include groups \cite{macmel}
which are finite as well as groups which are
infinite, both  connected
and not$-$connected.
The finite (compact) subgroups of $\; {\rm S}{\rm
O}(2)\times {\rm S}{\rm O}(2) \;$ were constructed explicitly in \cite{mel}.
The infinite not$-$connected compact
subgroups  of  $\; {\rm S}{\rm
O}(2)\times {\rm S}{\rm O}(2) \;$  were described  in detail in \cite{macmel}.
Here, we study in more detail, the infinite connected
(compact) subgroups of $\; {\rm S}{\rm O}(2)\times {\rm S}{\rm
O}(2) \;$. We start by explicitly constructing all (i.e. compact
and not$-$compact) infinite connected subgroups of $\; {\rm S}{\rm
O}(2)\times {\rm S}{\rm O}(2) \;$ and then we isolate the compact
infinite connected subgroups of $\; {\rm S}{\rm O}(2)\times {\rm
S}{\rm O}(2) . \;$ Firstly, we prove the following:

\vspace{0.5cm}

\begin{thrm}
\label{hyyyyhyhyhyhyiiokjklo} The connected subgroups of ${\rm
S}{\rm O}(2)\times {\rm S}{\rm O}(2)$ are classified as follows:
\begin{enumerate}
\item ${\rm S}{\rm O}(2)\times {\rm S}{\rm O}(2)$ \ is the only
two$-$dimensional connected subgroup.

\item The one$-$dimensional connected subgroups are: \be {\rm (a)}
\quad  \left(  \left( \begin{array}{cc}
\cos({\rm q}_{\rm o}\theta) & \sin({\rm q}_{\rm o}\theta) \\
\!\!\!\!\!-\sin({\rm q}_{\rm o}\theta) & \cos({\rm q}_{\rm
o}\theta)
\end{array}
\right),
 \left( \begin{array}{cc}
\cos({\rm p}_{\rm o}\theta) & \sin({\rm p}_{\rm o}\theta) \\
\!\!\!\!\!-\sin({\rm p}_{\rm o}\theta) & \cos({\rm p}_{\rm
o}\theta)
\end{array}
\right) \right) \equiv {\rm S}^{1}_{( {\rm q}_{\rm o}, {\rm
p}_{\rm o})}, \ee \be {\rm (b)} \quad \left(  \left(
\begin{array}{cc}
\cos\theta & \sin\theta \\
\!\!\!\!\!-\sin\theta & \cos\theta
\end{array}
\right),
 \left( \begin{array}{cc}
1 & 0 \\
0 & 1
\end{array}
\right) \right),
 \left( \left( \begin{array}{cc}
1 & 0 \\
0 & 1
\end{array}
\right), \left( \begin{array}{cc}
\cos\theta & \sin\theta \\
\!\!\!\!\!-\sin\theta & \cos\theta
\end{array}
\right) \right), \ee

\noindent where $\theta \in R$ \ and \ ${{\rm p}_{\rm o}>0,{\rm
q}_{\rm o}}>0$ \ or \ ${\rm p}_{\rm o}>0,{\rm q}_{\rm o}<0$. The
pair $({\rm q}_{\rm o},{\rm p}_{\rm o})$ is unique for each
subgroup. \item The identity element is the only zero$-$dimensional
connected subgroup.

\end{enumerate}
\end{thrm}
\begin{pf}
Connected subgroups are uniquely determined by their Lie algebras.
This is the content of a theorem given in Helgason (1962) \cite{Helgason}, which
states: (i) if H is a Lie subgroup of $\mathcal G$ \ its Lie
algebra dH is subalgebra of the Lie algebra d${\mathcal G}$ of
$\mathcal G$, \ and (ii) each Lie subalgebra of d${\mathcal G}$ is
the Lie algebra of exactly one connected subgroup of ${\mathcal
G}$.

\noindent
A basis for the Lie algebra of ${\rm S}{\rm O}(2)\times {\rm
S}{\rm O}(2)$ is given by the generators \be g_{1}= \left(  \left(
\begin{array}{cc}
0 & 1\\
\!\!\!\!-1 & 0
\end{array}
\right),
 \left( \begin{array}{cc}
0 & 0 \\
0 & 0
\end{array}
\right) \right),\;\;\quad
 g_{2}=\left( \left( \begin{array}{cc}
0 & 0 \\
0 & 0
\end{array}
\right), \left( \begin{array}{cc}
0 & 1\\
\!\!\!\!-1 & 0
\end{array}
\right) \right). \ee Therefore, there is only one two$-$dimensional
connected subgroup of ${\rm S}{\rm O}(2)\times {\rm S}{\rm O}(2)$,
\ which is \ ${\rm S}{\rm O}(2)\times {\rm S}{\rm O}(2)$ itself:
\begin{eqnarray}
{\rm e}^{\theta g_{1}}\cdot {\rm e}^{\phi g_{2}} & = &
 \left(  \left( \begin{array}{cc}
\cos\theta & \sin\theta \\
\!\!\!\!-\sin\theta & \cos\theta
\end{array}
\right),
 \left( \begin{array}{cc}
1 & 0 \\
0 & 1
\end{array}
\right) \right)\cdot \left(  \left( \begin{array}{cc}
1 & 0\\
0 & 1
\end{array}
\right),
 \left( \begin{array}{cc}
\cos\phi & \sin\phi\\
\!\!\!\!\!\!-\sin\phi & \cos\phi
\end{array}
\right) \right) \nonumber \\
 & & \left(  \left( \begin{array}{cc}
\cos\theta & \sin\theta \\
\!\!\!\!-\sin\theta & \cos\theta
\end{array}
\right),
 \left( \begin{array}{cc}
\cos\phi & \sin\phi \\
\!\!\!\!\!-\sin\phi & \cos\phi
\end{array}
\right) \right).
\end{eqnarray}
A basis for each one$-$dimensional Lie subalgebra of the Lie algebra
of ${\rm S}{\rm O}(2)\times {\rm S}{\rm O}(2)$ is given by the
generator: \be \label{sky} {\rm g}={\rm q}g_{1}+{\rm p}g_{2}=
\left(  \left( \begin{array}{cc}
0 & {\rm q} \\
\!\!\!\!\!-{\rm q} & 0
\end{array}
\right),
 \left( \begin{array}{cc}
0 & {\rm p} \\
\!\!\!\!\!-{\rm p} & 0
\end{array}
\right) \right), \ee where q, p are real numbers so  that $({\rm
q},{\rm p}) \neq (0,0).$ All one$-$dimensional connected subgroups
of ${\rm S}{\rm O}(2)\times {\rm S}{\rm O}(2)$ are obtained by
exponentiation: \be \label{sinepho} {\rm e}^{\theta {\rm g}}=
\left(  \left( \begin{array}{cc}
\cos({\rm q}\theta) & \sin({\rm q}\theta) \\
\!\!\!\!\!-\sin({\rm q}\theta) & \cos({\rm q}\theta)
\end{array}
\right),
 \left( \begin{array}{cc}
\cos({\rm p}\theta) & \sin({\rm p}\theta) \\
\!\!\!\!\!-\sin({\rm p}\theta) & \cos({\rm p}\theta)
\end{array}
\right) \right). \ee

({\rm a}) \quad   We shall consider first the case ${\rm q}{\rm
p}\neq 0$.\ Different bases of the same Lie subalgebra produce
the same group. In particular, (q,p) and ($-$q,$-$p) originate  the
same subgroup. The same applies to (q,$-$p) and ($-$q,p). Therefore,
only ${\rm p}>0,{\rm q}>0$ \  and ${\rm p}>0,
 {\rm q}<0$ have  to be considered. We now choose a representative
basis vector for each Lie subalgebra. For our later convenience,
when p/q is rational, we choose p,q to be the relatively prime
numbers ${\rm p}_{\rm o},{\rm q}_{\rm o}$, \ which satisfy ${\rm
p}_{\rm o}/{\rm q}_{\rm o}={\rm p}/{\rm q}$. When ${\rm p}/{\rm
q}=\gamma$  is irrational, we can consider  for example $({\rm
q}_{\rm o},{\rm p}_{\rm o})$  to be the point of interesection of
the open half$-$line from the origin $\delta {\rm g}$, $\delta \in
R^{+}$, with the half circle ${\rm q}^{2}+{\rm p}^{2}=1,{\rm
p}>0$. This choice having been made, $({\rm q}_{\rm o},{\rm
p}_{\rm o})=(1/\sqrt{1+\gamma^{2}}, |\gamma|/\sqrt{1+\gamma^{2}})$
for ${\rm p}>0, \ {\rm q}>0$ and $({\rm q}_{\rm o},{\rm p}_{\rm
o})=(-1/\sqrt{1+\gamma^{2}}, |\gamma|\sqrt{1+\gamma^{2}})$ for
${\rm p}>0, \ {\rm q}<0$.

({\rm b}) \quad The two remaining cases, ${\rm p}=0,{\rm q}\neq 0$
and ${\rm q}=0,{\rm p}\neq 0$, produce two more one-dimensional
subgroups, $({\rm SO(2)},{\rm I})$ and $({\rm I},{\rm SO(2)})$
respectively,where I is the $2\times 2$ unit  matrix.

Finally, for ${\rm p}={\rm q}=0$,\  Eqs. (\ref{sky}) and
(\ref{sinepho}) give the identity, as the only zero$-$dimensional
connected subgroup. This completes the proof.
\end{pf}
The one$-$dimensional subgroups are bounded. In order to find those
which are closed and therefore compact, firstly we prove the
following
\begin{prop}
The map h  defined by \be h:(R\times R)/(Z\times Z)\longrightarrow
{\rm S}{\rm O}(2)\times {\rm S}{\rm O}(2), \ee \be h((\frac
{\theta}{2\pi},\frac {\phi}{2\pi})+ Z\times Z)= \left(  \left(
\begin{array}{cc}
\cos\theta & \sin\theta \\
\!\!\!\!-\sin\theta & \cos\theta
\end{array}
\right),
 \left( \begin{array}{cc}
\cos\phi & \sin\phi \\
\!\!\!\!\!-\sin\phi & \cos\phi
\end{array}
\right) \right) \ee where $\theta, \phi \in R$ and  $Z$ denotes
the set of integers, is a homeomorphism.The quotient set $(R\times
R)/(Z\times Z) $ is endowed with the quotient topology and ${\rm
S}{\rm O}(2)\times {\rm S}{\rm O}(2)$ is equipped with the usual
topology.
\end{prop}
\begin{pf}
We define the maps $p$ and $\chi$ as follows:


\begin{figure}[ht]
\begin{center}
\begin{picture}(160,60)
\put(30,50){$R\times R$} \put(65,53){\vector(1,0){30}}
\put(76,58){$\chi$} \put(100,50){${\rm S}{\rm O}(2)\times {\rm
S}{\rm O}(2)$} \put(46,47){\vector(1,-1){30}} \put(50,5){$(R\times
R)/(Z\times Z)$} \put(86,17){\vector(1,1){30}} \put(53,27){$p$}
\put(103,24){$h$}
\end{picture}
\caption{\label{figure1} ${\rm S}{\rm O}(2)\times {\rm S}{\rm
O}(2)$ is homeomorphic to $(R\times R)/(Z\times Z).$}
\end{center}
\end{figure}
\vspace*{-11mm} \be \chi(\frac {\theta}{2\pi},\frac {\phi}{2\pi})=
\left(  \left( \begin{array}{cc}
\cos\theta & \sin\theta \\
\!\!\!\!-\sin\theta & \cos\theta
\end{array}
\right),
 \left( \begin{array}{cc}
\cos\phi & \sin\phi \\
\!\!\!\!\!-\sin\phi & \cos\phi
\end{array}
\right) \right), \ee \be p(\frac {\theta}{2\pi},\frac
{\phi}{2\pi})= (\frac {\theta}{2\pi},\frac {\phi}{2\pi})+Z\times
Z, \ee where $\theta,\phi \in R$. Since $\chi$ is a function , it
is clear that $[(\frac {\theta}{2\pi},\frac
{\phi}{2\pi})]=\{\chi^{-1}(g,h)|(g,h)\in SO(2)\times SO(2)\}$ is a
partition of $R\times R$. If $(\frac {\theta}{2\pi},\frac
{\phi}{2\pi}) \in R \times R$, $p((\frac {\theta}{2\pi},\frac
{\phi}{2\pi}))= [(\frac {\theta}{2\pi},\frac {\phi}{2\pi})]=
\chi^{-1}(\chi(\frac {\theta}{2\pi},^{}\frac {\phi}{2\pi})).$
 The map $\chi$ is evidently onto. The set
${\rm S}{\rm O}(2)\times {\rm S}{\rm O}(2)$, which is a subset of
$SL(2,R)\times SL(2,R) \subset R^{8}$, where $R^{8}$ is equipped
with the standard metric toplogy, inherits a relative topology. If
$R\times R$ is endowed with the usual topology the map $\chi$ is
continuous. Indeed, if $\chi$ is considered to be a function from
$R^{2}$ to $R^{8}$, it is continuous since the component functions
of $\chi$ are continuous. It is not difficult to show that the map
$\chi$ remains continuous, even if its range is restricted to the
subset ${\rm S}{\rm O}(2)\times {\rm S}{\rm O}(2)$, when it is
equipped with the relative topology.
If $(R\times R)/(Z\times Z)$ is endowed with the quotient
topology, then the well defined map $h$ is continuous
(\cite{aspects}, p.204). Furthermore, the map $h$ is one$-$to$-$one
and onto. Since $(R\times R)/(Z\times Z)$ is compact in the
quotient topology and ${\rm S}{\rm O}(2)\times {\rm S}{\rm O}(2)$
is Hausdorff in the relative topology, the map $h$ is a
homeomorphism (\cite{kelley}, p.141) and the proposition is
proved.
\end{pf}
\begin{prop}
The one$-$dimensional connected subgroups of \ ${\rm S}{\rm
O}(2)\times {\rm S}{\rm O}(2)$ with ${\rm q_{o}}{\rm p_{o}}\neq
~0$,\ ${\rm q_{o}}/{\rm p_{o}}$ irrational, are not compact in the
relative topology.
\end{prop}
\begin{pf}
We consider the subset of $R\times R$ \be \label{gata} {\rm
Q}=\left\{ (\frac {\theta {\rm q}}{2\pi} + {\rm t}, \frac {\theta
{\rm p}}{2\pi} + {\rm k})|\theta \in R, \quad {\rm t},{\rm k} \in
Z \right\} \ee where {\rm q},\ {\rm p} are fixed real numbers, so
that ${\rm q}{\rm p }\neq ~0$  and  ${\rm q}/{\rm p}$ is
irrational. For any $({\rm x},{\rm y}) \in R \times R$  we now
consider the distance \be {\rm d}=\left | \left | ({\rm x},{\rm
y})-(\frac {\theta {\rm q}}{2\pi} + {\rm t}, \frac {\theta {\rm
p}}{2\pi} + {\rm k}) \right | \right|. \ee We  can always write in
a non$-$unique way: $\theta=\frac {2 \pi}{ {\rm p}}(\phi +{\rm
j})|\phi \in R,{\rm j} \in Z$, and express {\rm d} accordingly \be
\label{xamogelo} {\rm d}=\left | \left | ( {\rm x}- \frac {\rm
q}{\rm p} \phi - \frac {\rm q}{\rm p} {\rm j}- {\rm t}, {\rm
y}-\phi-{\rm j}- {\rm k}) \right | \right |. \ee In expression
(\ref {xamogelo}), {\rm x}, {\rm y}, {\rm p}, {\rm q} are fixed,
while $\phi$, j, t, k can vary. The choice $~{\phi}= {\rm y}, {\rm
j}={\rm k}$ makes  the second component vanish \be {\rm d}=\left |
\left | ( {\rm x}- \frac {\rm q}{\rm p} {\rm y} - \frac {\rm
q}{\rm p} {\rm k}- {\rm t}, 0) \right | \right |=\left  | {\rm x}-
\frac {\rm q}{\rm p} {\rm y} - \frac {\rm q}{\rm p} {\rm k}- {\rm
t} \right |. \ee The subgroup $Z+ \frac{\rm q}{\rm p}Z $  of  $R$
is dense in $R$, if and only if  $\frac{\rm q}{\rm p}$ is
irrational \cite{hardy}. Henceforth, for any $\epsilon > 0$
there always exist integers k, t, so that ${\rm d}< \epsilon$.
Consequently, the set Q is dense in $R \times R$. The projection
map $p$ (see figure 1.1) is continuous and onto. Therefore, the
image $p({\rm Q})$ is dense in $(R\times R)/(Z \times Z)$. Since
$p({\rm Q})$ is  a proper subset of $(R\times R)/(Z\times Z)$,
$p({\rm Q})$ is not closed in the quotient topology. We conclude
that $h(p({\rm Q}))$ is not a closed and consequently not a
compact subgroup of ${\rm S}{\rm O}(2)\times {\rm S}{\rm O}(2)$ in
the relative topology. If $h(p({\rm Q}))$ were a closed subgroup
of ${\rm S}{\rm O}(2)\times {\rm S}{\rm O}(2)$, $h^{-1}(h(p({\rm
Q}))=p({\rm Q})$ would be a closed subgroup of $(R\times
R)/(Z\times Z)$, since $h^{-1}$, which is a homeomorphism, is a
closed map. This completes the proof.
\end{pf}
\begin{prop}
The one$-$dimensional connected subgroups of \ ${\rm S}{\rm
O}(2)\times {\rm S}{\rm O}(2)$ with ${\rm q_{o}}{\rm p_{o}}\neq
~0$,\ ${\rm q_{o}}/{\rm p_{o}}$ rational and either ${\rm
q_{o}}\neq 0, \ {\rm p_{o}}=0$  or ${\rm p_{o}}\neq 0, \ {\rm q_{o}}=0$
are  compact in the relative topology.
\end{prop}
\begin{pf}The elements of the set Q, which is defined by (Eq.(\ref{gata})),
span straight lines on the cartesian plane $R\times R$. These
straight lines, in standard cartesian coordinates x, y, if  ${\rm
q}{\rm p}\neq ~0$, are given by \be \label{kalo} {\rm y}=\frac
{{\rm p}}{{\rm q}} {\rm x} - \frac {{\rm p}}{{\rm q}} {\rm t}
+{\rm k}. \ee Consequently,  they are parallel with slope p/q.
When ${\rm q}\neq 0, \ {\rm p}=0$  and  ${\rm p}\neq 0, \ {\rm q}=0$,
they are given respectively by \be \label{kalo1} {\rm y}= {\rm k}
\quad {\rm and} \quad {\rm x}={\rm t}, \ee where ${\rm k},{\rm t}
\in Z $. Since the straight lines (\ref{kalo}) and (\ref{kalo1})
are labelled by a countable set, the set Q with either ${\rm
q}{\rm p}\neq ~0$ or ${\rm q}\neq 0,{\rm p}=0$, \ \ ${\rm p}\neq
0,{\rm q}=0$ is a proper subset of $R\times R$. Firstly, we
consider  case

(a) \ ${\rm q}{\rm p}\neq ~0$,\ ${\rm q}/{\rm p}$ rational. The
ordinates of the intersections of the straight lines (\ref{kalo})
with the y$-$axis are given by \be {\rm y}_{({\rm t},{\rm k})}= -
\frac{{\rm p}}{{\rm q}} {\rm t} +{\rm k}. \ee For $ \frac {\rm
p}{\rm q}$ rational, the subgroup $Z+ \frac  {\rm q}{\rm p} Z$ of
$R$, is not dense in $R$ and consequently  is a discrete subgroup
of $R$. Therefore, the distance between two adjacent straight
lines (\ref{kalo}) cannot become arbitrarily small for any choice
of the integers k and t. We conclude that the complement of Q is
open. Henceforth, Q itself is closed. Q=$p^{-1}(p(\rm Q))$ is
closed and therefore $p(\rm Q)$ is a closed subset
(\cite{aspects},p.204) of $(R\times R)/(Z\times Z) $ in the
quotient topology. Since $h$ is closed, $h(p(\rm Q))$ is a closed
and therefore a compact subgroup of ${\rm S}{\rm  O}(2)\times {\rm
S}{\rm  O}(2)$ in the relative topology. Now we consider case

(b) \  ${\rm q}\neq 0,{\rm p}=0$  or ${\rm p}\neq 0,{\rm q}=0$.
With a proof similar to the one  given in case (a) we can show
that $h(p(\rm Q))$ is a closed and therefore a compact subgroup of
${\rm S}{\rm  O}(2)\times {\rm S}{\rm  O}(2)$ in the relative
topology. This completes the proof.
\end {pf}

The infinite not$-$connected compact
subgroups  of  $\; {\rm S}{\rm
O}(2)\times {\rm S}{\rm O}(2) \;$  were described  in detail in \cite{macmel}.
In particular the following Theorem was proved:

\begin{thrm}
The infinite not$-$connected subgroups of $ \ {\rm S}{\rm O}(2)\times {\rm S}{\rm
O}(2)$ are all one$-$dimensional and they are  classified as
follows:
\begin{enumerate}

\item  \begin{eqnarray}  \left(  \left( \begin{array}{cc}
\cos({\rm q}_{\rm o}\theta) & \sin({\rm q}_{\rm o}\theta) \\
\!\!\!\!\!-\sin({\rm q}_{\rm o}\theta) & \cos({\rm q}_{\rm
o}\theta)
\end{array}
\right),
 \left( \begin{array}{cc}
\cos({\rm p}_{\rm o}\theta + \frac {2 \pi}{N}i) & \sin({\rm p}_{\rm o}\theta + \frac {2 \pi}{N}i) \\
\!\!\!\!\!-\sin({\rm p}_{\rm o}\theta + \frac {2 \pi}{N}i) &
\cos({\rm p}_{\rm o}\theta + \frac {2 \pi}{N}i) \end{array}
\right) \right) &\equiv& \nonumber \\  H(N,{\rm q}_{\rm o},{\rm
p}_{\rm o})\simeq C_{N} \times {\rm S}^{1}_{( {\rm q}_{\rm o},
{\rm p}_{\rm o})}, \quad  &&
\end{eqnarray}

\item \be
 \left( \left( \begin{array}{cc}
\cos( \frac{2 \pi}{N}i) & \sin( \frac{2 \pi}{N}i) \\
\!\!\!\!\!\!-\sin ( \frac{2 \pi}{N}i) & \cos( \frac{2 \pi}{N}i)
\end{array}
\right), \left( \begin{array}{cc}
\cos\theta & \sin\theta \\
\!\!\!\!\!-\sin\theta & \cos\theta
\end{array}
\right) \right) \simeq C_{N} \times S^{1}, \quad  \ee

\item \be
 \left(  \left(
\begin{array}{cc}
\cos\theta & \sin\theta \\
\!\!\!\!\!-\sin\theta & \cos\theta
\end{array}
\right),
 \left( \begin{array}{cc}
\cos( \frac{2 \pi}{N}i)  & \sin( \frac{2 \pi}{N}i) \\
\!\!\!\!\!\!-\sin ( \frac{2 \pi}{N}i)  & \cos( \frac{2 \pi}{N}i)
\end{array}
\right) \right) \simeq S^{1} \times C_{N}, \quad  \ee \noindent

\end{enumerate}
where, the positive  integer $\;N\;$ satisfies  $\; N> 1$, $
\; i \in \{0,1,2,...,N-1 \}$, $\theta \in R,$ \ and where,\
the coprime integers ${{\rm p}_{\rm o}, \ {\rm q}_{\rm o}}$ satisfy
${{\rm p}_{\rm o}>0,{\rm q}_{\rm o}}>0$ \ or \ ${\rm p}_{\rm
o}>0,{\rm q}_{\rm o}<0$.
\end{thrm}

\section{Infinite  Potential Little Groups}

\label{s2}

In order to carry out step (II), we must find, for each infinite
subgroup $\;S \subseteq {\rm S}{\rm O}(2) \times {\rm S}{\rm
O}(2)\;$ the space $\;A(S)\;$ of its invariant functions.
We start by proving the following:
\begin{prop}
\label{kiloikjuhbbnbvfg} There are no non$-$zero functions \ $ \phi_
{\rm o} \in {L^{2}}^{'}(\mathcal P,\lambda ,R)$ \ fixed under {\rm :}

\begin{enumerate}

\item ${\rm S}{\rm O}(2) \times {\rm S}{\rm O}(2)$.

\item $ \; S^{1} \times C_{N} \; $ where $ \; N \geq 1 \;$.

\item $ \; C_{N} \times S^{1} \;$ where $ \; N \geq 1 \;$.

\item $H(N, {\rm q}_{\rm o}, {\rm p}_{\rm o})$ where one of the
numbers  $ \; N,$ $ {\rm q}_{\rm o}$, $ {\rm p}_{\rm
o} $ is even.
\end{enumerate}
\end{prop}
\begin{pf}

\quad 1. \quad $ \left(  \left( \begin{array}{cc}
\cos\theta & \sin\theta \\
\!\!\!\!-\sin\theta & \cos\theta
\end{array}
\right),
 \left( \begin{array}{cc}
\cos\varphi & \sin\varphi \\
\!\!\!\!\!-\sin\varphi & \cos\varphi
\end{array}
\right) \right) $

\vspace*{3mm}

We require that  $\phi_ {\rm o}(x,y)$  be fixed under
 ${\rm S}{\rm O}(2) \times {\rm S}{\rm O}(2)$.  Taking into account
(\ref{aliki}), (\ref{donkey}), (\ref{mule}), and (\ref{tzatziki}),
the condition on $\phi_ {\rm o}(x,y)$ reads
$$
\!\!\!\!\!\!\!\!\!\!\!\!\!\!\!\!\!\!\!\!\!\!\!\!\!\!\!\!\!\!\!\!\!\!\!\!\!\!\!
\!\!\!\!\!\!\!\! \phi_ {\rm o}(x,y) =  (T'(g,h)\phi_ {\rm o})(x,y)
=k_{g}^{-3}(x)s_{g}(x)k_{h}^{-3}(y)s_{h}(y)\phi(xg,yh)
$$
\be \label{piperia} \;\;\;\;\;\;\;\;\;\;\;\;\;\;\;\;\; =\frac {x
\sin \theta + \cos \theta}{|x \sin \theta + \cos \theta|} \frac {y
\sin \varphi  + \cos \varphi}{|y \sin \varphi  + \cos \varphi|}
\phi_ {\rm o} \left ( \frac {x \cos \theta - \sin \theta} {x \sin
\theta + \cos \theta}, \frac {y \cos \varphi - \sin \varphi} {y
\sin \varphi  + \cos \varphi} \right ) \ee where $- \infty <
\theta, \varphi < + \infty $. As an alternative co$-$ordinatization
of $P_{1}(R) \times P_{1}(R)$ we can use angular co$-$ordinates $
\rho,  \sigma $ related to the projective co$-$ordinates $x, y$ by
the equations \be \label{kolokithi} x= \cot \frac{\rho}{2} \qquad
y=\cot \frac {\sigma}{2}, \quad - \infty < \rho, \sigma < + \infty
. \ee Hereafter, \ $ \varsigma (x) \equiv s (\sin x) $.\ By using
(\ref{kolokithi}) the invariant condition (\ref{piperia}) on
$\phi_{\rm o}(x,y)$  becomes \be \label{choma} \phi_{\rm o}(\rho +
2 \theta, \sigma + 2 \varphi)= \varsigma \left ( \frac {\rho + 2
\theta}{2} \right ) \varsigma \left ( \frac {\rho}{2} \right )
\varsigma \left ( \frac {\sigma + 2 \varphi}{2} \right ) \varsigma
\left ( \frac {\sigma}{2} \right ) \phi_{\rm o}(\rho, \sigma), \ee
where $\phi_{\rm o}(\rho, \sigma) \stackrel{\rm def}{=}\phi_{\rm
o}(x,y)| _{x= \cot \frac{\rho}{2},y=\cot \frac {\sigma}{2}}$. We
must further require that $\phi_{\rm o}(\rho, \sigma)$ be doubly
periodic \be \label{xorto} \phi_{\rm o}(\rho + 2 \pi,
\sigma)=\phi_{\rm o}(\rho, \sigma), \qquad \phi_{\rm o}(\rho,
\sigma + 2 \pi)=\phi_{\rm o}(\rho, \sigma). \ee Consequently, \be
\label{troxos} \phi_{\rm o}(\rho + 2\theta +2\pi, \sigma + 2
\varphi)= \phi_{\rm o}(\rho + 2\theta , \sigma + 2 \varphi). \ee
However, from (\ref {choma}) we have \be \label{karotsa} \phi_{\rm
o}(\rho + 2(\theta +\pi), \sigma + 2 \varphi)= - \phi_{\rm o}(\rho
+ 2\theta , \sigma + 2 \varphi). \ee From (\ref{troxos}) and
(\ref{karotsa}) we conclude: \be \phi_{\rm o}(\rho + 2\theta ,
\sigma + 2 \varphi)=0. \ee Hence, \ $\phi_{\rm o}$ \ vanishes as
an element of ${L^{2}}^{'}(\mathcal P,\lambda ,R)$.

\vspace*{4mm}

\hspace*{21mm}2. \quad $
 \left(  \left(
\begin{array}{cc}
\cos\theta & \sin\theta \\
\!\!\!\!\!-\sin\theta & \cos\theta
\end{array}
\right),
 \left( \begin{array}{cc}
\cos( \frac{2 \pi}{N}i)  & \sin( \frac{2 \pi}{N}i) \\
\!\!\!\!\!\!-\sin ( \frac{2 \pi}{N}i)  & \cos( \frac{2 \pi}{N}i)
\end{array}
\right) \right) $




\vspace*{4mm} \noindent The condition that \ $\phi_{\rm o}$ \ be
fixed under $\;  S^{1} \times C_{N} \;$
\ is \begin{eqnarray} \phi_ {\rm o}(x,y) & = & \frac {x \sin
\theta + \cos \theta}{|x \sin \theta + \cos \theta|} \frac {y \sin
( \frac{2 \pi}{N}i) + \cos ( \frac{2 \pi}{N}i)}{|y \sin ( \frac{2
\pi}{N}i) + \cos ( \frac{2 \pi}{N}i)|} \nonumber \\
&&\phi_ {\rm o} \left ( \frac {x \cos \theta - \sin \theta} {x
\sin \theta + \cos \theta}, \frac {y \cos ( \frac{2 \pi}{N}i) -
\sin ( \frac{2 \pi}{N}i)} {y \sin ( \frac{2 \pi}{N}i) + \cos (
\frac{2 \pi}{N}i)} \right ).
\end{eqnarray} Using (\ref{kolokithi}) we find: \be
\label{pirinas} \phi_{\rm o}(\rho + 2 \theta, \sigma +2( \frac{2
\pi}{N}i))= \varsigma \left ( \frac {\rho + 2 \theta}{2} \right )
\varsigma \left ( \frac {\rho}{2} \right ) \varsigma \left ( \frac
{\sigma + 2 ( \frac{2 \pi}{N}i)}{2} \right ) \varsigma \left (
\frac {\sigma}{2} \right ) \phi_{\rm o}(\rho, \sigma). \ee From
(\ref{pirinas}) we have: \be \phi_{\rm o}(\rho + 2(\theta +\pi),
\sigma +2( \frac{2 \pi}{N}i) )= - \phi_{\rm o}(\rho + 2\theta ,
\sigma +2( \frac{2 \pi}{N}i)). \ee The periodicity condition
(\ref{xorto}) then implies: \be \phi_{\rm o}(\rho + 2\theta ,
\sigma +2( \frac{2 \pi}{N}i))=0, \ee where \ $ -\infty <
\rho, \theta , \sigma < + \infty $. Henceforth, we find that \
$\phi_{\rm o}$ \ vanishes as an element of ${L^{2}}^{'}(\mathcal P,\lambda ,R)$.

\vspace*{4mm}

\hspace*{21mm}3. \quad $
 \left( \left( \begin{array}{cc}
\cos( \frac{2 \pi}{N}i) & \sin( \frac{2 \pi}{N}i) \\
\!\!\!\!\!\!-\sin ( \frac{2 \pi}{N}i) & \cos( \frac{2 \pi}{N}i)
\end{array}
\right), \left( \begin{array}{cc}
\cos\theta & \sin\theta \\
\!\!\!\!\!-\sin\theta & \cos\theta
\end{array}
\right) \right) \simeq C_{N} \times S^{1}, $


\vspace*{4mm}

This case is very similar to the previous one. The condition that
$\phi_{\rm o}$  be fixed under  $ \; C_{N} \times S^{1} \;$ \ is:
\be \label{kordoni} \phi_{\rm o}(\rho + 2(\frac{2 \pi}{N}i),
\sigma + 2 \theta)= \varsigma \left ( \frac {\rho + 2(
\frac{2\pi}{N}i)}{2} \right ) \varsigma \left ( \frac {\rho}{2}
\right ) \varsigma \left ( \frac {\sigma + 2 \theta}{2} \right )
\varsigma \left ( \frac {\sigma}{2} \right ) \phi_{\rm o}(\rho,
\sigma). \ee Eq. (\ref{kordoni}) yields \ $ \phi_{\rm o}(\rho + 2(
\frac{2\pi}{N}i) , \sigma + 2(\theta +\pi))= - \phi_{\rm o}(\rho +
2( \frac{2\pi}{N}i) , \sigma + 2\theta)$ \, and, using \
(\ref{xorto}), we conclude as before that \ $\phi_{\rm o}$ \
vanishes as an element of ${L^{2}}^{'}(\mathcal P,\lambda ,R)$.

\vspace*{4mm}

\hspace*{21mm}4.  \quad $  \left(  \left(
\begin{array}{cc}
\cos({\rm q}_{\rm o}\theta) & \sin({\rm q}_{\rm o}\theta) \\
\!\!\!\!\!-\sin({\rm q}_{\rm o}\theta) & \cos({\rm q}_{\rm
o}\theta)
\end{array}
\right),
 \left( \begin{array}{cc}
\cos({\rm p}_{\rm o}\theta + \frac {2 \pi}{N}i) & \sin({\rm p}_{\rm o}\theta + \frac {2 \pi}{N}i) \\
\!\!\!\!\!-\sin({\rm p}_{\rm o}\theta + \frac {2 \pi}{N}i) &
\cos({\rm p}_{\rm o}\theta + \frac {2 \pi}{N}i) \end{array}
\right) \right), $


\vspace*{4mm}

\noindent where one of the numbers \ $ \; N,$ $ \; {\rm
q}_{\rm o}, $  $ {\rm p}_{\rm o} $  is even . An element
\ $\phi_{\rm o} \in A(\N)$ \ is invariant under $ \; H(N,{\rm
q}_{\rm o},{\rm p}_{\rm o})\simeq C_{N} \times {\rm S}^{1}_{( {\rm
q}_{\rm o}, {\rm p}_{\rm o})} \;$
\ if and only if \begin{eqnarray} \label{kidoni} \phi_{\rm o}(\rho
+ 2 {\rm q}_{\rm o}\theta, \sigma + 2{\rm p}_{\rm o} \theta + 2
(\frac{2 \pi}{N}i)) & = & \varsigma \left ( \frac {\rho + 2 {\rm
q}_{\rm o}\theta }{2} \right ) \varsigma \left ( \frac {\rho}{2}
\right ) \nonumber \\
&& \varsigma \! \left ( \frac {\sigma + 2{\rm p}_{\rm o} \theta +
2 (\frac{2 \pi}{N}i)}{2} \right ) \varsigma \! \left ( \frac
{\sigma}{2} \right ) \phi_{\rm o}(\rho, \sigma).
\end{eqnarray}
We distinguish two cases:

\begin{enumerate}
\item One of the numbers $\; {\rm p}_{{\rm o}}$, $\;{\rm
q}_{{\rm o}} \;$ is even.

\vspace{0.4cm}

 From the periodicity condition (\ref{xorto}), we have: \be
\label{xarti} \phi_{\rm o}(\rho + 2{\rm q}_{\rm o}(\theta + \pi),
\sigma + 2{\rm p}_{\rm o}(\theta + \pi)+ 2 (\frac{2 \pi}{N}i))=
\phi_{\rm o}(\rho + 2{\rm q}_{\rm o}\theta, \sigma + 2{\rm p}_{\rm
o}\theta+ 2 (\frac{2 \pi}{N}i)). \ee Since one of the numbers \
${\rm q}_{\rm o}, {\rm p}_{\rm o}$ \ is even and the other is odd,
Eq. (\ref{kidoni}) gives: \be \label{eukaliptos} \phi_{\rm o}(\rho
+ 2{\rm q}_{\rm o}(\theta + \pi), \sigma + 2{\rm p}_{\rm o}(\theta
+ \pi)+ 2 (\frac{2 \pi}{N}i))= - \phi_{\rm o}(\rho + 2{\rm q}_{\rm
o}\theta, \sigma + 2{\rm p}_{\rm o}\theta+ 2 (\frac{2 \pi}{N}i)).
\ee By combining (\ref{xarti}) and (\ref{eukaliptos}) we find again
that \ $\phi_{\rm o}$ \ vanishes as an element of ${L^{2}}^{'}(\mathcal P,\lambda ,R)$.

\item The number $ \; N \;$ is even.

\vspace{0.4cm}

 From the periodicity condition (\ref{xorto}), we have:
 \be  \label{nmjkiiiolplo} \phi_{\rm o}(\rho + 2 {\rm q}_{\rm o}\theta, \sigma + 2{\rm
p}_{\rm o} \theta + 2 (\frac{2 \pi}{N}(i+\frac{N}{2})))= \phi_{\rm
o}(\rho + 2 {\rm q}_{\rm o}\theta, \sigma + 2{\rm p}_{\rm o}
\theta + 2 (\frac{2 \pi}{N}i)) \ . \ee Eq. (\ref{kidoni}) gives:
\be \label{kijuuuhhhyyhjukinbes}\phi_{\rm o}(\rho + 2 {\rm q}_{\rm
o}\theta, \sigma + 2{\rm p}_{\rm o} \theta + 2 (\frac{2
\pi}{N}(i+\frac{N}{2})))= - \phi_{\rm o}(\rho + 2 {\rm q}_{\rm
o}\theta, \sigma + 2{\rm p}_{\rm o} \theta + 2 (\frac{2 \pi}{N}i))
\ . \ee Combining (\ref{nmjkiiiolplo}) and
(\ref{kijuuuhhhyyhjukinbes}), we find again that \ $\phi_{\rm o}$
\ vanishes as an element of ${L^{2}}^{'}(\mathcal P,\lambda ,R)$.
\end{enumerate}
This completes  the proof .
\end{pf}


\noindent
Two remarks are now in order.

 \noindent Firstly, we recall
 the notion of fundamental region.
  Let $\Omega_{\phi_{{\rm o}}}$ \ be the set of periods of \
$ {\phi_{{\rm o}}} $. $ \Omega_{\phi_{{\rm o}}} $ \ forms a group
isomorphic to $ Z \times Z $ \ under the usual addition. We may
regard \ $\Omega_{\phi_{{\rm o}}}$ \ as acting on $ R^{2} $ as a
transformation group, each $(2 \pi m, 2 \pi n) \in
\Omega_{\phi_{{\rm o}}}$, where \ $ (m,n) \in Z \times Z $, \
inducing the translation \be {\rm t}_ {(2 \pi m, 2 \pi n)} \ : \ (
\rho , \sigma ) \longmapsto ( \rho + 2 \pi m, \sigma + 2 \pi n)
\ee of $ R^{2}$. A closed, connected subset  P of \ $ R^{2}$ \ is
defined to be a fundamental region for \ $\Omega_{\phi_{{\rm o}}}$
\ if

(i) for each \ $ (\rho, \sigma ) \in R^{2} $, P contains at least
one point in the same \ $\Omega_{\phi_{{\rm o}}}-$ orbit as
\hspace*{1.09cm} $ (\rho, \sigma )$

(ii) no two points in the interior of P are in the same
$\Omega_{\phi_{{\rm o}}}$$-$ orbit.


\noindent
Secondly,
it is an immediate
corollary of section 6 in \cite{macmel}
that two groups
$\;H(N_{1}, {\rm q}_{1},{\rm p}_{1}) \simeq   C_{N_{1}}\times{\rm
S}^{1}_{({\rm q}_{1},{\rm p}_{1})}$ and $\;H(N_{2}, {\rm
q}_{2},{\rm p}_{2})\simeq   C_{N_{2}}\times{\rm S}^{1}_{({\rm
q}_{2},{\rm p}_{2})}\;$ are equal if and only if $\;
C_{N_{1}}=C_{N_{2}}\;$ and $\;{\rm S}^{1}_{({\rm q}_{1},{\rm
p}_{1})}={\rm S}^{1}_{({\rm q}_{2},{\rm p}_{2})}.\;$ Taking into
account Theorem \ref{hyyyyhyhyhyhyiiokjklo} and Proposition
\ref{kiloikjuhbbnbvfg}
we conclude that there is no loss of
generality if we assume that $N$ is odd and that
${\rm q}, \ {\rm p}$ in $\;H(N,{\rm q},{\rm p})\;$
are odd and relatively prime. With this assumption
the equality $\;H(N_{1}, {\rm
q}_{1},{\rm p}_{1}) = H(N_{2}, {\rm q}_{2},{\rm p}_{2})\;$ implies
$\; {\rm q}_{1}={\rm q}_{2} \;$ and $\; {\rm p}_{1}={\rm
p}_{2}.\;$ Hereafter we assume that $\; N \;$ is odd and
$\;{\rm q},{\rm p}\;$ are odd and relatively prime.

\noindent
We
show now that there are non$-$zero
functions $\phi_{\rm o} \in {L^{2}}^{'}(\mathcal P,\lambda ,R)$
 fixed under each one of the remaining compact infinite subgroups
$\; H (N,{\rm q},{\rm p}) \;$
of \ ${\rm SO(2) \times \rm SO(2)}$.
Taking into account the results of Theorem
\ref{hyyyyhyhyhyhyiiokjklo} we
distinguish four cases. We also find
elementary regions for the action
of the groups $\;H(N,{\rm q},{\rm
p})$ on the torus $ {\rm T} \simeq P_{1}(R)\times P_{1}(R)$.

\begin{thrm}
\label{iiiuolokjuuuihyjuik} Every  subgroup \ $\; H (N,{\rm
q}_{{\rm o}},{\rm p}_{{\rm o}}) \;$  of \ ${\rm SO(2) \times \rm
SO(2)}$, where $ \; N \;$ is odd and $ \;{\rm q}_{{\rm o}},{\rm
p}_{{\rm o}} \;$ are odd and relatively prime,
 has non$-$zero fixed functions \ $\phi_{\rm o} \in
A(\N)$
The fixed functions and the corresponding little groups are given
in TABLE 1. The spaces of fixed functions \ $\vartheta _{{\rm q}N'}
(u_{i})$, where i=1,2,3,4 and, \ $ u_{1}=\sigma - {\rm q} \rho $,
$ u_{2}=\rho  - {\rm q} \sigma $,\ $ u_{3}=\sigma + {\rm q}\rho
$,\ $ u_{4}=-\rho - {\rm q}\sigma $,\ defined on arcs \ ${\rm
A}_{i{\rm q}_{{\rm o}}N'}$ on the torus \ $ P_{1}(R)\times
P_{1}(R)$, \ are the Hilbert spaces \ ${\rm L}^{2}_{i{\rm
q}N'}({\rm A}_{i{\rm q}_{{\rm o}}N'}, \lambda, R)$ \ of square-
integrable functions with respect to the Lebesgue measure \ $
\lambda $ \ along the arcs \ ${\rm A}_{i{\rm q}_{{\rm o}}N'}$. The
arcs  $ \; {\rm A}_{i{\rm q}_{{\rm o}}N'} \;$ are given by
$\;v_{i}=0 \ , \ 0 \leq u_{i} \leq \frac{2\pi}{{\rm q}_{{\rm
o}}N'}\;$, $\;i=1,2,3,4.\;$ The coordinates $\;v_{i}\;$ are
defined as follows: $\; v_{1}={\rm q}\sigma + \rho, \;$ $\;
v_{2}={\rm q}\rho + \sigma, \;$ $\; v_{3}={\rm q}\sigma - \rho,
\;$ $\; v_{4}=\sigma - {\rm q}\rho . \;$ When $\; i=1 \;$ and $\;
i=3 ,\;$ we define  $\;N=N'l\;$ and $\;{\rm q}_{{\rm o}}={\rm
q}_{{\rm o}}'l,\;$ where $\;l={\rm gcd}(N,{\rm q}_{{\rm o}}).\;$
When $\; i=2 \;$ and $\; i=4 ,\;$ we define $\;N=N'l\;$ and
$\;{\rm p}_{{\rm o}}={\rm p}_{{\rm o}}'l,\;$ where $\;l={\rm
gcd}(N,{\rm p}_{{\rm o}}).\;$

\center {\rm TABLE  \ \ 1}

\hspace*{0.8 in}{\rm fixed functions} $\phi_{\rm o}$ \hspace*{1.1
in}{\rm corresponding little groups} $L_{\phi_{\rm o}}$


\begin{enumerate}
\item $\phi_{{\rm o}1}(\rho,\sigma)= \varsigma \left ( \frac {-
\rho}{2} \right ) \varsigma \left ( \frac {\sigma}{2} \right )
\vartheta_{{\rm q}N'} (\sigma - {\rm q} \rho)$, ${\rm q} \equiv
\frac {{\rm p}_{\rm o}}{{\rm q}_{\rm o}}$ \hspace*{0.2in}$H( N,
{\rm q}_{\rm o}, {\rm p}_{\rm o})$, ${\rm p}_{\rm o}>{\rm q}_{\rm
o}>0 \ \ {\rm or} \ \ {\rm p}_{\rm o}= \hspace*{3.11in}  {\rm
q}_{\rm o}=1.$

\item $\phi_{{\rm o}2}(\rho,\sigma)= \varsigma \left ( \frac {-
\sigma}{2} \right ) \varsigma \left ( \frac {\rho}{2} \right )
\vartheta_{{\rm q}N'} (\rho  - {\rm q} \sigma)$, ${\rm q} \equiv
\frac {{\rm p}_{\rm o}}{{\rm q}_{\rm o}}$ \hspace*{0.2in}$H(N ,
{\rm p}_{\rm o}, {\rm q}_{\rm o})$, ${\rm p}_{\rm o}>{\rm q}_{\rm
o}>0 \ \ {\rm or} \ \ {\rm p}_{\rm o}= \hspace*{3.11in}{\rm
q}_{\rm o}=1.$

\item $\phi_{{\rm o}3}(\rho,\sigma)= \varsigma \left ( \frac {
\sigma}{2} \right )\varsigma \left ( \frac {\rho}{2} \right )
\vartheta_{{\rm q}N'} (\sigma + {\rm q} \rho)$  ,
 $ {\rm q} \equiv \frac {{\rm p}_{\rm o}}{{\rm q}_{\rm o}}$
\hspace*{0.2in}$H(N,- {\rm q}_{\rm o}, {\rm p}_{\rm o})$, ${\rm
p}_{\rm o}>{\rm q}_{\rm o}>0 \ \  {\rm or} \ \ {\rm p}_{\rm
o}=\hspace*{3.06in} {\rm q}_{\rm o}=1.$

\item $\phi_{{\rm o}4}(\rho,\sigma)= \varsigma \left ( \frac {
\sigma}{2} \right ) \varsigma \left ( \frac {\rho}{2} \right )
\vartheta_{{\rm q}N'} (-\rho - {\rm q} \sigma)$  ,
 $ {\rm q} \equiv \frac {{\rm p}_{\rm o}}{{\rm q}_{\rm o}}$
\hspace*{0.2in}$\!\!H(N,- {\rm p}_{\rm o}, {\rm q}_{\rm o})$,
${\rm p}_{\rm o}>{\rm q}_{\rm o}>0 \ \  {\rm or} \ \ {\rm p}_{\rm
o}= \hspace*{3.06in}{\rm q}_{\rm o}=1.$

\end{enumerate}
\end{thrm}

\begin{pf}
The action of  \ ${\rm S}^{1}_{( {\rm q}_{\rm o}, {\rm p}_{\rm
o})}$ \ on the torus $ {\rm T} \simeq P_{1}(R)\times P_{1}(R)$
partitions \ $ {\rm T}^{2}$ \ into a union of disjoint closed
orbits. Each of these orbits turns around the torus \ ${\rm q}_{\rm
o}$ \ times along its meridians, clockwise for \ ${\rm q}_{\rm o}$
\ positive and anti-clockwise for \ ${\rm q}_{\rm o}$ \ negative,
and goes \ ${\rm p}_{\rm o}$ \ times along its latitudes,
clockwise (anti-clockwise) for \ ${\rm p}_{\rm o}$ \ positive
(negative). An orbit of the action of $\; H(N, {\rm q}_{\rm o},
{\rm p}_{\rm o})\;$ on the torus $ \;  {\rm T} \simeq P_{1}(R)\times P_{1}(R)\; $ consists of
$ \; N \;$ orbits of   the action of $\;{\rm S}^{1}_{( {\rm
q}_{\rm o}, {\rm p}_{\rm o})}$ \ on the torus $ {\rm T} \simeq P_{1}(R)\times P_{1}(R). \;$

\vspace*{4mm}

\hspace*{1mm}1.  \ We look for functions \ $ \phi_{{\rm o}1}$ \
which are fixed under \ $H(N,{\rm q}_{\rm o}, {\rm p}_{\rm o})$,\
where ${\rm p}_{\rm o}>{\rm q}_{\rm o}>~0 $\ \ {\rm or} \ \ $~{\rm
p}_{\rm o}= {\rm q}_{\rm o}=1.$ \ The invariance condition \
$(T'(g,h)\phi_{{\rm o}1})(\rho,\sigma)=\phi_{{\rm
o}1}(\rho,\sigma)$ gives:
\begin{eqnarray}
\label{hjjjddddsnh} \phi_{{\rm o}1}(\rho + 2 {\rm q}_{\rm
o}\theta, \sigma + 2{\rm p}_{\rm o} \theta + 2(\frac{2\pi}{N}i)) &
= & \varsigma \left ( \frac {\rho + 2 {\rm q}_{\rm o}\theta }{2}
\right ) \varsigma \left ( \frac {\rho}{2} \right ) \varsigma
\left ( \frac {\sigma + 2{\rm p}_{\rm o} \theta +
2(\frac{2\pi}{N}i)}{2}  \right ) \nonumber \\
&& \label{korena} \varsigma \left ( \frac {\sigma}{2} \right )
\phi_{{\rm o}1}(\rho, \sigma). \end{eqnarray}



\noindent Among  the infinite fundamental regions for \
$\Omega_{\phi_{{\rm o}1}}$,  the usual  representative is the
square region $ ~  \{ ( \rho , \sigma ) \in R^{2}: 0 \leq \rho ,
\sigma \leq 2\pi \}.$ In the problem under consideration, a
fundamental region which is spanned by the $\;{\rm S}^{1}({\rm
q}_{{\rm o}},{\rm p}_{{\rm o}})-{\rm orbits}\;$ is a more natural
choice.
Such a fundamental region is the region $\;{\rm F}_{1}\;$ which is
depicted in Figure \ref {melas}. We can easily verify that this
region satisfies (i) and (ii). Moreover, although a fundamental
region for $\Omega_{\phi_{{\rm o}1}}$, is not unique, its area is
unique and we can readily find that the area of the region we
chose is equal to  $4 \pi^{2} $. We now make the following change
of co$-$ordinates: \be \label{kratiras} u_{1}=\sigma - {\rm q} \rho
\qquad v_{1}={\rm q} \sigma + \rho, \quad {\rm q}\equiv \frac
{{\rm p}_{\rm o}}{{\rm q}_{\rm o}}. \ee The new coordinate $u_{1}$
remains constant along each \ ${\rm S}^{1}_{( {\rm q}_{\rm o},
{\rm p}_{\rm o})}$$-$orbit. The lines of the two families of lines
$u_{1}$=constant and $v_{1}$=constant are straight and
perpendicular to each other. In the $ u_{1},v_{1} $ coordinates
the fundamental region we chose is given by \be 0 \leq u_{1} \leq
\frac {2 \pi}{{\rm q}_{\rm o}}, \qquad 0 \leq v_{1} \leq  {2
\pi}{\rm q}_{\rm o}(1+{\rm q}^{2}). \ee We define \be
\label{karagouna} \psi_{{\rm o}1}(u_{1},v_{1})=\phi_{{\rm o}1}
\left ( \frac {v_{1}-{\rm q}u_{1}}{1+{\rm q}^{2}}, \frac
{v_{1}{\rm q}+u_{1}} {1+{\rm q}^{2}} \right ) \equiv \phi_{{\rm
o}1}(\rho, \sigma). \ee The function $ \; \psi_{{\rm
o}1}(u_{1},v_{1}) \;$ satisfies a periodicity condition. Recall
that the function $\; \phi_{\rm o}(\rho , \sigma) \;$ is doubly
periodic (Eq. (\ref{xorto})):
$$  \phi_{\rm o}(\rho + 2 \pi, \sigma)=\phi_{\rm
o}(\rho, \sigma), \quad  \quad \phi_{\rm o}(\rho, \sigma + 2
\pi)=\phi_{\rm o}(\rho, \sigma). $$ The last periodicity
conditions are equivalent to the following one: \be
\label{mikklloooolkiolkio} \phi_{\rm o}(\rho + 2 m\pi, \sigma +
2n\pi)=\phi_{\rm o}(\rho, \sigma) \ , \ee where $ \; m \ {\rm and}
\ n \;$ are integers. Eq. (\ref{xorto}) imply: \be
\label{huhhhuuuhu} \psi_{{\rm o}1}(u_{1} - 2 \pi{\rm q},
v_{1}+2\pi)=\psi_{{\rm o}1}(u_{1}, v_{1}), \quad  \quad \psi_{{\rm
o}1}(u_{1}+2\pi, v_{1} + 2 \pi{\rm q})=\psi_{{\rm o}1}(u_{1},
v_{1}).  \ee Eq. (\ref{mikklloooolkiolkio}) implies: \be
\label{koiiuiuyhyjuik} \psi_{{\rm o}1}(u_{1}+a \frac{2\pi}{{\rm
q}_{{\rm o}}} \ , \ v_{1}+ b \frac{2\pi}{{\rm q}_{{\rm o}}}) =
\psi_{{\rm o}1},  \ee  where $\; a \;$ and $\; b \;$ are
integers. The periodicity conditions (\ref{huhhhuuuhu}) and
(\ref{koiiuiuyhyjuik}) are equivalent. It turns out that on some
occasions is more convenient to use one rather than the other.

The functional Eq (\ref{hjjjddddsnh}) in the new co$-$ordinates
$u_{1},v_{1}$ is written as follows:
\begin{eqnarray}
\psi_{{\rm o}1}(u_{1}+ \frac{4 \pi}{N}\nu ,v_{1}+(1+{\rm
q}^{2})2{\rm q}_{\rm o}\theta+ \frac{4\pi}{N}{\rm q}\nu) & = &
\varsigma \left ( \frac {v_{1}-{\rm q}u_{1} + 2 {\rm q}_{\rm o}
\theta (1 + {\rm q}^{2}) } { 2 (1 + {\rm q}^{2})} \right )
\varsigma \left ( \frac {v_{1}-{\rm q}u_{1}}{2(1+{\rm q}^{2})}
\right )  \nonumber \\ &  & \varsigma \left ( \frac {v_{1} {\rm q}
+ u_{1} + 2{\rm p}_ {\rm o}  \theta (1 + {\rm q}^{2})+ \frac{4
\pi}{N}\nu(1+{\rm q}^{2}) }  { 2 (1 + {\rm q}^{2})} \right )
\nonumber \\ && \label{kollokkijuhy} \varsigma \left ( \frac
 {v_{1} {\rm q}+u_{1}}
{2(1+{\rm q}^{2})} \right ) \psi_{{\rm o}1}(u_{1},v_{1}).
\end{eqnarray}
When \ $ v_{1}=0 $ and $ \nu =0, $   the last equation yields:
\begin{eqnarray} \psi_{{\rm o}1}(u_{1},(1+{\rm q}^{2})2{\rm q}_{\rm
o}\theta) & = & \varsigma \left ( \frac {-{\rm q}u_{1} + 2 {\rm
q}_{\rm o} \theta (1 + {\rm q}^{2}) } { 2 (1 + {\rm q}^{2})}
\right ) \varsigma \left ( \frac {-{\rm q}u_{1} } { 2 (1 + {\rm
q}^{2})}
\right ) \nonumber \\
&& \varsigma \left ( \frac { u_{1} + 2{\rm p}_ {\rm o}  \theta (1
+ {\rm q}^{2}) } { 2 (1 + {\rm q}^{2})} \right ) \varsigma \left (
\frac { u_{1} } { 2 (1 + {\rm q}^{2})} \right )
 \psi_{{\rm o}1}(u_{1},0).
\end{eqnarray}
The last Eq. reads
\begin{eqnarray} \psi_{{\rm o}1}(u_{1},v_{1}) & = &
\varsigma \left ( \frac {{\rm q}u_{1} -v_{1} } { 2 (1 + {\rm
q}^{2})} \right ) \varsigma \left ( \frac {{\rm q}u_{1} } { 2 (1 +
{\rm q}^{2})}
\right ) \nonumber \\
&& \label{kiooioioiojiklo} \varsigma \left ( \frac { u_{1} + {\rm
q}v_{1} } { 2 (1 + {\rm q}^{2})} \right ) \varsigma \left ( \frac
{ u_{1} } { 2 (1 + {\rm q}^{2})} \right )
 \psi_{{\rm o}1}(u_{1},0).
\end{eqnarray}
Using the periodicity condition (\ref{koiiuiuyhyjuik}) one without
loss of generality can restrict $ \; u_{1} \;$ to the range of
values $ \; u_{1} \in [0,\frac{2\pi}{{\rm q}_{{\rm o}}}] .\;$ In
this case, Eq. (\ref{kiooioioiojiklo}) gives
\be \label{mosxos} \psi_{{\rm o}1}(u_{1},v_{1})= \varsigma \left (
\frac { {\rm q}u_{1} - v_{1}  }{2(1+{\rm q}^{2})} \right )
\varsigma \left ( \frac
 { u_{1} + {\rm q}v_{1} }
{2(1+{\rm q}^{2})} \right )  \psi_{{\rm o}1}(u_{1},0) \ . \ee We
find now which periodicity condition the function $ \; \psi_{{\rm
o}1}(u_{1},0) \;$ satisfies. Imposing on Eq. (\ref{mosxos}) the
periodicity condition $ \;
 \psi_{{\rm o}1}(u_{1} - 2 \pi{\rm q}, v_{1}+2\pi)=\psi_{{\rm o}1}(u_{1}, v_{1})
 \;$ (Eq.(\ref{huhhhuuuhu})) we obtain
 \be
 \label{hujkiiiloplo}
 \psi_{{\rm o}1}(u_{1}-2\pi{\rm q},0)=-\psi_{{\rm o}1}(u_{1},0).
 \ee
Imposing on Eq. (\ref{mosxos}) the periodicity condition $ \;
\psi_{{\rm o}1}(u_{1}+2\pi, v_{1} + 2 \pi{\rm q})=\psi_{{\rm
o}1}(u_{1}, v_{1}) \ \; $ (Eq.(\ref{huhhhuuuhu})) we obtain \be
\label{lolooiiuyhy} \psi_{{\rm o}1}(u_{1}+2\pi,0)=-\psi_{{\rm
o}1}(u_{1},0).
\ee Combining the periodicity conditions (\ref{hujkiiiloplo}) and
(\ref{lolooiiuyhy}) we obtain \be
\label{qqqqqqqqqjhdaaaaa}\psi_{{\rm o}1}(u_{1}+ \frac{2 \pi}{{\rm
q}_{{\rm o}}},0)=-\psi_{{\rm o}1}(u_{1},0).  \ee The periodicity
condition (\ref{qqqqqqqqqjhdaaaaa}) is equivalent to the
periodicity conditions (\ref{hujkiiiloplo}) and
(\ref{lolooiiuyhy}). From the periodicity condition
(\ref{qqqqqqqqqjhdaaaaa}) one can easily derive the periodicity
conditions (\ref{hujkiiiloplo}) and (\ref{lolooiiuyhy}). We
explain now how one from  the periodicity conditions
(\ref{hujkiiiloplo}) and (\ref{lolooiiuyhy}) can derive the
periodicity condition (\ref{qqqqqqqqqjhdaaaaa}). From Eqs.
(\ref{hujkiiiloplo}) and (\ref{lolooiiuyhy}) we obtain
\be\label{yyyyyhgdvbawe}\psi_{{\rm o}1}(u_{1}+{\rm m}2\pi{\rm q} +
{\rm n}2 \pi)= - \psi_{{\rm o}1}(u_{1}), \quad  \ee where one of
the integers m,n is even and the other is odd. Since for given
relatively prime numbers ${\rm p}_{\rm o}, \ {\rm q}_{\rm o}$ there
always exist integers m, n $-$ one even and one odd $-$ for which $
{\rm m} \frac {{\rm p}_{\rm o}}{{\rm q}_{\rm o}} + {\rm n}= \frac
{1} {{\rm q}_{\rm o}} $, Eq. (\ref{yyyyyhgdvbawe}) gives
$$
\psi_{{\rm o}1}(u_{1}+ \frac{2 \pi}{{\rm q}_{{\rm
o}}},0)=-\psi_{{\rm o}1}(u_{1},0), \  $$ which is Eq.
(\ref{qqqqqqqqqjhdaaaaa}).
The function $ \; \psi_{{\rm o}1}(u_{1},v_{1}) \;$ (Eq.
(\ref{mosxos})) satisfies Eq. (\ref{kollokkijuhy}) and
substituting the expression (\ref{mosxos}) into the Eq.
(\ref{kollokkijuhy}) we obtain
 \be
 \label{mjuuuuiuiuiuiuiu}
 \psi_{{\rm o}1}(u_{1}+ \frac{4 \pi}{N} \nu ,0)=\psi_{{\rm
o}1}(u_{1},0).
 \ee
Combining the periodicity conditions (\ref{qqqqqqqqqjhdaaaaa}) and
(\ref{mjuuuuiuiuiuiuiu}) we obtain \be \psi_{{\rm
o}1}(u_{1}+m\frac{2\pi}{{\rm q}_{{\rm
o}}}+\frac{4\pi}{N}\nu,0)=-\psi_{{\rm o}1}(u_{1},0),  \ee where
$\; m \;$ is an {\it odd} integer. We note that $ \; m \frac{2
\pi}{{\rm q}_{{\rm o}}}+ \frac{4\pi}{N}\nu=\frac{2\pi}{{\rm
q}_{{\rm o}}N}(mN+{\rm q}_{{\rm o}}(2\nu)). \;$ Let $\;l={\rm
gcd}(N,{\rm q}_{{\rm o}}) \;$ be the greatest common divisor of
$\; N \;$ and $\; {\rm q}_{{\rm o}}. \;$
We have $\;N=N'l\;$ and $\;{\rm q}_{{\rm o}}={\rm q}_{{\rm
o}}'l\;$, where $\;N'\;$ and $\;{\rm q}_{{\rm o}}'\;$ are coprime.
Taking into account the last Eqs. we obtain:  $ \; m \frac{2
\pi}{{\rm q}_{{\rm o}}}+ \frac{4\pi}{N}\nu=\frac{2\pi}{{\rm
q}_{{\rm o}}N}l(mN'+{\rm q}_{{\rm o}}'(2\nu)). \;$ Since $\; N'
\;$ and $\; {\rm q}_{{\rm o}}' \;$ are coprime, there always
exist two integers, one odd $\;(m)\;$ and one even $\;(2\nu)\;$
such that $\;mN'+{\rm q}_{{\rm o}}'(2\nu)=1.\;$ We conclude that
the periodicity condition (\ref{mjuuuuiuiuiuiuiu}) reads
\be\label{mikiloooopoijijik} \psi_{{\rm
o}1}(u_{1}+\frac{2\pi}{{\rm q}_{{\rm o}}N}l,0)=-\psi_{{\rm
o}1}(u_{1},0),  \ee or equivalently, \be \psi_{{\rm
o}1}(u_{1}+\frac{2\pi}{{\rm q}_{{\rm o}}N'},0)=-\psi_{{\rm
o}1}(u_{1},0).  \ee


From Eqs. (\ref {kratiras}), (\ref {karagouna}) and (\ref {mosxos})
we conclude that the solution to the functional Eq.(\ref {korena})
is given almost everywhere by: \be \label{zorbas} \phi_{{\rm
o}1}(\rho,\sigma)= \varsigma \left ( \frac {- \rho}{2} \right )
\varsigma \left ( \frac {\sigma}{2} \right ) \vartheta_{{\rm q}N'}
(\sigma - {\rm q} \rho), \ee where $\vartheta_{{\rm q}N'} (u_{1})
\equiv \psi_{{\rm o}1}(u_{1},0)$. For a given triad $\;(N,{\rm
q}_{{\rm o}},{\rm p}_{{\rm o}})\;$,
the functions \ $\phi_{{\rm o}1}(\rho,\sigma)$ \ given by
(\ref{zorbas}) form the maximal invariant closed vector space \
$A(H(N,{\rm q}_{\rm o}, {\rm p}_{\rm o})) \in A(\N)$ of \ $H(N,
{\rm q}_{\rm o}, {\rm p}_{\rm o})$, where either ${\rm p}_{\rm
o}>{\rm q}_{\rm o}>0 \quad {\rm or} \quad {\rm p}_{\rm o}={\rm
q}_{\rm o}=1.$ Note that since the numbers $\; {\rm q}_{{\rm o}}
\;$ and $\; {\rm p}_{{\rm o}} \;$ are chosen to be coprime,
specific choice of $\; {\rm q} \;$  is equivalent to determining
uniquely the numbers $\; {\rm q}_{{\rm o}} \;$ and $\; {\rm
p}_{{\rm o}}. $
\noindent The square$-$integrable condition on $\;\phi_{{\rm
o}1}(\rho,\sigma)\;$ gives:
\begin{eqnarray}
\int_{{\rm F}_{1}} \int | \phi_{{\rm o}1} (\rho,\sigma)|^{2}{\rm
d}\rho \, {\rm d}\sigma  & = & \int_{0}^{2\pi{\rm q}_{\rm o}(1+
{\rm q}^{2})} \int_{0}^{2\pi / {\rm q}_{\rm o}} |\psi_{{\rm
o}1}(u_{1},v_{1})|^{2} \; | \left | \frac { \partial (\rho,
\sigma)} { \partial (u_{1}, v_{1})} \right | | \;
{\rm d}u_{1} \, {\rm d}v_{1} \nonumber \\
& = &  2\pi {\rm q}_{\rm o} \int _{0}^{2\pi / {\rm q}_{\rm o}} |
\vartheta_{{\rm q}N'} (u_{1}) |^ {2} {\rm d}u_{1} \nonumber \\
& = &  2\pi {\rm q}_{\rm o}N' \int _{0}^{2\pi / {\rm q}_{\rm o}N'}
| \vartheta_{{\rm q}N'} (u_{1}) |^ {2} {\rm d}u_{1} < + \infty,
\end{eqnarray}
where ${\rm F}_{1}$ denotes the fundamental region depicted in
Fig. \ref{melas} and $| \left | \frac { \partial (\rho, \sigma)}
{\partial (u_{1}, v_{1})} \right | | $ is the absolute value of
the Jacobian of the co$-$ordinate transformation $ (u_{1},v_{1})
\longmapsto (\rho, \sigma)$. Therefore, $A(H(N, {\rm q}_{\rm o},
{\rm p}_{\rm o}))$, as a Hilbert space, is isomorphic to ${\rm
L}^{2}_{1{\rm q}N'}({\rm A}_{1{\rm q}_{{\rm o}}N'}, \lambda, R)$,
the Hilbert space of the real valued functions $ \vartheta _{{\rm
q}N'}(u_{1})$ which are defined on the arc ${\rm A}_{1{\rm
q}_{{\rm o}}N'}$ $(v_{1}=0, \; 0 \leq u_{1} \leq \frac {2
\pi}{{\rm q}_{\rm o}N'})$ and are square$-$integrable with respect
to the Lebesgue measure $\; \lambda \;$ along this arc.
\vspace*{4mm}

\noindent We give now the examination of the remaining cases. At
many places,in each case, the mathematical treatment is similar to
the first case. We omit derivations which are similar to those
given in the first case.

\vspace{0.5cm}

\hspace*{1mm}2.
 \ We look for functions \ $ \phi_{{\rm o}2}$ \
which are fixed under \ $H(N,{\rm p}_{\rm o}, {\rm q}_{\rm o})$,\
where ${\rm p}_{\rm o}>{\rm q}_{\rm o}>~0 $\ \ {\rm or} \ \ $~{\rm
p}_{\rm o}= {\rm q}_{\rm o}=1.$ \ The invariance condition \
$(T'(g,h)\phi_{{\rm o}2})(\rho,\sigma)=\phi_{{\rm
o}2}(\rho,\sigma)$ gives:
\begin{eqnarray}
 \phi_{{\rm o}2}(\rho + 2 {\rm p}_{\rm o}\theta, \sigma + 2{\rm
q}_{\rm o} \theta + 2(\frac{2\pi}{N}i)) & = & \varsigma \left (
\frac {\rho + 2 {\rm p}_{\rm o}\theta }{2} \right ) \varsigma
\left ( \frac {\rho}{2} \right ) \varsigma \left ( \frac {\sigma +
2{\rm q}_{\rm o} \theta +
2(\frac{2\pi}{N}i)}{2}  \right ) \nonumber \\
&& \label{jhhdsxcbnhe} \varsigma \left ( \frac {\sigma}{2} \right
) \phi_{{\rm o}2}(\rho, \sigma). \end{eqnarray}



\noindent
A fundamental region for $\Omega_{\phi_{{\rm o}2}}$ which is
spanned by the  $\;{\rm S}^{1}({\rm p}_{{\rm o}},{\rm q}_{{\rm
o}})-{\rm orbits}\;$ is the region $\;{\rm F}_{2}\;$ which is
depicted in Figure \ref {melas}. The fundamental regions $\;{\rm
F}_{1}\;$ and $\;{\rm F}_{2}\;$ are symmetric with respect to the
axis $\;\rho=\sigma.\;$
We now make the following change of co$-$ordinates: \be
\label{kiiiioouaasre} u_{2}=\rho - {\rm q} \sigma, \qquad
v_{2}={\rm q} \rho + \sigma, \quad \rm{where} \quad  {\rm q}\equiv \frac {{\rm
p}_{\rm o}}{{\rm q}_{\rm o}}. \ee The new coordinate $u_{2}$
remains constant along each \ ${\rm S}^{1}_{( {\rm p}_{\rm o},
{\rm q}_{\rm o})}$$-$orbit. The lines of the two families of lines
$u_{2}$=constant and $v_{2}$=constant are straight and
perpendicular to each other. In the $ u_{2}, \ v_{2} $ coordinates
the fundamental region we chose is given by \be 0 \leq u_{2} \leq
\frac {2 \pi}{{\rm q}_{\rm o}}, \qquad 0 \leq v_{2} \leq  {2
\pi}{\rm q}_{\rm o}(1+{\rm q}^{2}). \ee We define \be
\label{ooooooasaew} \psi_{{\rm o}2}(u_{2},v_{2})=\phi_{{\rm o}2}
\left ( \frac {u_{2}+{\rm q}v_{2}}{1+{\rm q}^{2}}, \frac
{v_{2}-{\rm q}u_{2}} {1+{\rm q}^{2}} \right ) \equiv \phi_{{\rm
o}2}(\rho, \sigma). \ee The function $ \; \psi_{{\rm
o}2}(u_{2},v_{2}) \;$ satisfies the following periodicity
condition:   \be \label{kilokiiiijuhynbv} \psi_{{\rm o}2}(u_{2}+a
\frac{2\pi}{{\rm q}_{{\rm o}}} \ , \ v_{2}+ b \frac{2\pi}{{\rm
q}_{{\rm o}}}) = \psi_{{\rm o}2}(u_{2},v_{2}), \ee where $\; a
\;$ and $\; b \;$ are integers. This last periodicity condition
is equivalent to the following periodicity conditions: \be
\label{shhhhyujnmbdf} \psi_{{\rm o}2}(u_{2} - 2 \pi{\rm q},
v_{2}+2\pi)=\psi_{{\rm o}2}(u_{2}, v_{2}),  \quad \psi_{{\rm
o}2}(u_{2}+2\pi, v_{2} + 2 \pi{\rm q})=\psi_{{\rm o}2}(u_{2},
v_{2}).  \ee


\noindent The functional Eq (\ref {jhhdsxcbnhe}) in the new
co$-$ordinates $u_{2},v_{2}$ is written as follows:
\begin{eqnarray}
\psi_{{\rm o}2}(u_{2}-{\rm q} \frac{4 \pi}{N}\nu ,v_{2}+(1+{\rm
q}^{2})2{\rm q}_{\rm o}\theta+ \frac{4\pi}{N}\nu) & = & \varsigma
\left ( \frac {u_{2}+{\rm q}v_{2} + 2 {\rm p}_{\rm o} \theta (1 +
{\rm q}^{2}) } { 2 (1 + {\rm q}^{2})} \right ) \varsigma \left (
\frac {u_{2}+{\rm q}v_{2}}{2(1+{\rm q}^{2})} \right )  \nonumber
\\ &  & \varsigma \left ( \frac {v_{2}- {\rm q} u_{2} + 2{\rm q}_
{\rm o}  \theta (1 + {\rm q}^{2})+ \frac{4 \pi}{N}\nu(1+{\rm
q}^{2}) }  { 2 (1 + {\rm q}^{2})} \right ) \nonumber \\ &&
\label{njkmloiuuuy} \varsigma \left ( \frac
 {v_{2} - {\rm q}u_{2}}
{2(1+{\rm q}^{2})} \right ) \psi_{{\rm o}2}(u_{2},v_{2}).
\end{eqnarray}
By setting  $\; v_{2}=0, \;$  $\; \nu =0 \;$ and by
 restricting $ \; u_{2} \;$ to the range of
values $ \; u_{2} \in [0,\frac{2\pi}{{\rm q}_{{\rm o}}}], \;$ we
obtain from Eq. (\ref{njkmloiuuuy}) \be \label{ujjiikiooolpollmk}
\psi_{{\rm o}2}(u_{2},v_{2})= \varsigma \left ( \frac { u_{2}
+{\rm q} v_{2} }{2(1+{\rm q}^{2})} \right ) \varsigma \left (
\frac
 { {\rm q}u_{2} - v_{2} }
{2(1+{\rm q}^{2})} \right )  \psi_{{\rm o}2}(u_{2},0) \ . \ee The
function $ \; \psi_{{\rm o}2}(u_{2},0) \;$ satisfies  the
following periodicity condition: \be\label{jukilomnj} \psi_{{\rm
o}2}(u_{2}+\frac{2\pi}{{\rm q}_{{\rm o}}N}l,0)=-\psi_{{\rm
o}2}(u_{2},0), \ee where $\;l={\rm gcd}(N,{\rm p}_{{\rm
o}}).\;$ Eq. (\ref{jukilomnj}) is equivalent to:
\be\label{miolokiolkiol} \psi_{{\rm o}2}(u_{2}+\frac{2\pi}{{\rm
q}_{{\rm o}}N'},0)=-\psi_{{\rm o}2}(u_{2},0),  \ee where
$\;N=N'l,\;$ $\;{\rm p}_{{\rm o}}={\rm p}_{{\rm o}}'l,\;$ and
$\;N', \ {\rm p}_{{\rm o}}'\;$ are coprime.
From Eqs. (\ref {kiiiioouaasre}), (\ref {ooooooasaew}) and (\ref
{ujjiikiooolpollmk}) we conclude that the solution to the
functional Eq.(\ref {jhhdsxcbnhe}) is given almost everywhere by:
\be \label{uuuawwwwwsdcxdfvcd} \phi_{{\rm o}2}(\rho,\sigma)=
\varsigma \left ( \frac {-\sigma}{2} \right ) \varsigma \left (
\frac {\rho}{2} \right ) \vartheta_{{\rm q}N'} (\rho - {\rm q}
\sigma), \ee where $\vartheta_{{\rm q}N'} (u_{2}) \equiv
\psi_{{\rm o}2}(u_{2},0)$.
\noindent The square$-$integrable condition on $\;\phi_{{\rm
o}2}(\rho,\sigma)\;$ gives:
\be
\int _{0}^{2\pi / {\rm q}_{\rm o}N'} | \vartheta_{{\rm q}N'}
(u_{2}) |^ {2} {\rm d}u_{2} < + \infty, \ee
Therefore, $A(H(N, {\rm p}_{\rm o}, {\rm q}_{\rm o}))$, as a
Hilbert space, is isomorphic to ${\rm L}^{2}_{2{\rm q}N'}({\rm
A}_{2{\rm q}_{{\rm o}}N'}, \lambda, R)$, the Hilbert space of the
real valued functions $ \vartheta _{{\rm q}N'}(u_{2})$ which are
defined on the arc ${\rm A}_{2{\rm q}_{{\rm o}}N'}$ $(v_{2}=0, \;
0 \leq u_{2} \leq \frac {2 \pi}{{\rm q}_{\rm o}N'})$ and are
square$-$integrable with respect to the Lebesgue measure $\; \lambda
\;$ along this arc. \vspace*{4mm}

\vspace*{4mm}

\hspace*{1mm}3.  \

 \ We look for functions \ $ \phi_{{\rm o}3}$ \
which are fixed under \ $H(N,-{\rm q}_{\rm o}, {\rm p}_{\rm o})$,\
where ${\rm p}_{\rm o}>{\rm q}_{\rm o}>~0 $\ \ {\rm or} \ \ $~{\rm
p}_{\rm o}= {\rm q}_{\rm o}=1.$ \ The invariance condition \
$(T'(g,h)\phi_{{\rm o}3})(\rho,\sigma)=\phi_{{\rm
o}2}(\rho,\sigma)$ gives:
\begin{eqnarray}
 \phi_{{\rm o}3}(\rho - 2 {\rm q}_{\rm o}\theta, \sigma + 2{\rm
p}_{\rm o} \theta + 2(\frac{2\pi}{N}i)) & = & \varsigma \left (
\frac {\rho - 2 {\rm q}_{\rm o}\theta }{2} \right ) \varsigma
\left ( \frac {\rho}{2} \right ) \varsigma \left ( \frac {\sigma +
2{\rm p}_{\rm o} \theta +
2(\frac{2\pi}{N}i)}{2}  \right ) \nonumber \\
&& \label{hhgbbnvfcdxxx} \varsigma \left ( \frac {\sigma}{2}
\right ) \phi_{{\rm o}3}(\rho, \sigma). \end{eqnarray}

\noindent A fundamental region for $\Omega_{\phi_{{\rm o}3}}$
which is spanned by the  $\;{\rm S}^{1}(-{\rm q}_{{\rm o}},{\rm
p}_{{\rm o}})-{\rm orbits}\;$ is the region $\;{\rm F}_{3}\;$
which is depicted in Figure \ref {bolos}. The fundamental
regions $\;{\rm F}_{3}\;$ and $\;{\rm F}_{1}\;$ are symmetric with
respect to the axis $\;\rho=0.\;$
We now make the following change of co$-$ordinates: \be
\label{hjxssewdfrv} u_{3}=\sigma + {\rm q} \rho, \qquad v_{3}={\rm
q} \sigma - \rho, \quad \rm{where} \quad{\rm q}\equiv \frac {{\rm p}_{\rm o}}{{\rm
q}_{\rm o}}. \ee The new coordinate $u_{3}$ remains constant along
each \ ${\rm S}^{1}_{( -{\rm q}_{\rm o}, {\rm p}_{\rm o})}$$-$orbit.
The lines of the two families of lines $u_{3}$=constant and
$v_{3}$=constant are straight and perpendicular to each other. In
the $ u_{3},v_{3} $ coordinates the fundamental region we chose is
given by \be 0 \leq u_{3} \leq \frac {2 \pi}{{\rm q}_{\rm o}},
\qquad 0 \leq v_{3} \leq  {2 \pi}{\rm q}_{\rm o}(1+{\rm q}^{2}).
\ee We define \be \label{njbdddssserfgdhy} \psi_{{\rm
o}3}(u_{3},v_{3})=\phi_{{\rm o}3} \left ( \frac {{\rm
q}u_{3}-v_{3}}{1+{\rm q}^{2}}, \frac {u_{3}+{\rm q}v_{3}} {1+{\rm
q}^{2}} \right ) \equiv \phi_{{\rm o}3}(\rho, \sigma). \ee The
function $ \; \psi_{{\rm o}3}(u_{3},v_{3}) \;$ satisfies the
following periodicity condition:   \be \label{jknhdsnbmnkjjjjs}
\psi_{{\rm o}3}(u_{3}+a \frac{2\pi}{{\rm q}_{{\rm o}}} \ , \
v_{3}+ b \frac{2\pi}{{\rm q}_{{\rm o}}}) = \psi_{{\rm
o}3}(u_{3},v_{3}), \ee where $\; a \;$ and $\; b \;$ are
integers. This last periodicity condition is equivalent to the
following periodicity conditions: \be \label{njmbhfnnbdvgfhnn}
\psi_{{\rm o}3}(u_{3} + 2 \pi{\rm q}, v_{3}-2\pi)=\psi_{{\rm
o}3}(u_{3}, v_{3}) \quad , \quad \psi_{{\rm o}3}(u_{3}+2\pi, v_{3}
+ 2 \pi{\rm q})=\psi_{{\rm o}3}(u_{3}, v_{3}). \ee

\noindent The functional Eq (\ref {hhgbbnvfcdxxx}) in the new
co$-$ordinates $u_{3}, \ v_{3}$ is written as follows:
\begin{eqnarray}
\psi_{{\rm o}3}(u_{3}+ \frac{4 \pi}{N}\nu ,v_{3}+(1+{\rm
q}^{2})2{\rm q}_{\rm o}\theta+ {\rm q}\frac{4\pi}{N}\nu) & = &
\varsigma \left ( \frac {{\rm q}u_{3}-v_{3} + 2 {\rm q}_{\rm o}
\theta (1 + {\rm q}^{2}) } { 2 (1 + {\rm q}^{2})} \right )
\varsigma \left ( \frac {{\rm q}u_{3}-v_{3}}{2(1+{\rm q}^{2})}
\right )  \nonumber
\\ &  & \varsigma \left ( \frac {{\rm q}v_{3}+ u_{3} + 2{\rm p}_
{\rm o}  \theta (1 + {\rm q}^{2})+ \frac{4 \pi}{N}\nu(1+{\rm
q}^{2}) }  { 2 (1 + {\rm q}^{2})} \right ) \nonumber \\ &&
\label{kiiiiiijhffffdeer} \varsigma \left ( \frac
 {{\rm q}v_{3} +u_{3}}
{2(1+{\rm q}^{2})} \right ) \psi_{{\rm o}3}(u_{3},v_{3}).
\end{eqnarray}
By setting  $\; v_{3}=0, \;$  $\; \nu =0 \;$ and by
 restricting $ \; u_{3} \;$ to the range of
values $ \; u_{3} \in [0,\frac{2\pi}{{\rm q}_{{\rm o}}}], \;$ we
obtain from Eq. (\ref{kiiiiiijhffffdeer}) \be
\label{oiiyyrfhgyjtuknf} \psi_{{\rm o}3}(u_{3},v_{3})= \varsigma
\left ( \frac
 { {\rm q}u_{3} - v_{3} }
{2(1+{\rm q}^{2})} \right ) \varsigma \left ( \frac { u_{3} +{\rm
q} v_{3} }{2(1+{\rm q}^{2})} \right ) \psi_{{\rm o}3}(u_{3},0).
\ee The function $ \; \psi_{{\rm o}3}(u_{3},0) \;$ satisfies  the
following periodicity condition: \be\label{olkmjnhbyuiokl}
\psi_{{\rm o}3}(u_{3}+\frac{2\pi}{{\rm q}_{{\rm
o}}N}l,0)=-\psi_{{\rm o}3}(u_{3},0), \ee where $\;l={\rm
gcd}(N,{\rm q}_{{\rm o}}).\;$ Eq. (\ref{olkmjnhbyuiokl}) is
equivalent to: \be\label{oooirrrfdds} \psi_{{\rm
o}3}(u_{3}+\frac{2\pi}{{\rm q}_{{\rm o}}N'},0)=-\psi_{{\rm
o}3}(u_{3},0),  \ee where $\;N=N'l,\;$ $\;{\rm q}_{{\rm o}}={\rm
q}_{{\rm o}}'l,\;$ and $\;N', \ {\rm q}_{{\rm o}}'\;$ are coprime.

From Eqs. (\ref {hjxssewdfrv}), (\ref {njbdddssserfgdhy}) and (\ref
{oiiyyrfhgyjtuknf}) we conclude that the solution to the
functional Eq.(\ref {hhgbbnvfcdxxx}) is given almost everywhere
by: \be \label{oouybqqqreedaws} \phi_{{\rm o}3}(\rho,\sigma)=
\varsigma \left ( \frac { \rho}{2} \right ) \varsigma \left (
\frac {\sigma}{2} \right ) \vartheta_{{\rm q}N'} (\sigma + {\rm q}
\rho), \ee where $\vartheta_{{\rm q}N'} (u_{3}) \equiv \psi_{{\rm
o}3}(u_{3},0)$.
\normalfont \noindent The square$-$integrable condition on
$\;\phi_{{\rm o}3}(\rho,\sigma)\;$ gives: \be \int _{0}^{2\pi /
{\rm q}_{\rm o}N'} | \vartheta_{{\rm q}N'} (u_{3}) |^ {2} {\rm
d}u_{3} < + \infty, \ee Therefore, $A(H(N, - {\rm q}_{\rm o}, {\rm
p}_{\rm o}))$, as a Hilbert space, is isomorphic to ${\rm
L}^{2}_{3{\rm q}N'}({\rm A}_{3{\rm q}_{{\rm o}}N'}, \lambda, R)$,
the Hilbert space of the real valued functions $ \vartheta _{{\rm
q}N'}(u_{3})$ which are defined on the arc ${\rm A}_{3{\rm
q}_{{\rm o}}N'}$ $(v_{3}=0, \; 0 \leq u_{3} \leq \frac {2
\pi}{{\rm q}_{\rm o}N'})$ and are square$-$integrable with respect
to the Lebesgue measure $\; \lambda \;$ along this arc.
\vspace*{4mm}

\vspace*{4mm}

\hspace*{1mm}4.  \
 \ We look for functions \ $ \phi_{{\rm o}4}$ \
which are fixed under \ $H(N,-{\rm p}_{\rm o}, {\rm q}_{\rm o})$,\
where ${\rm p}_{\rm o}>{\rm q}_{\rm o}>~0 $\ \ {\rm or} \ \ $~{\rm
p}_{\rm o}= {\rm q}_{\rm o}=1.$ \ The invariance condition \
$(T'(g,h)\phi_{{\rm o}4})(\rho,\sigma)=\phi_{{\rm
o}4}(\rho,\sigma)$ gives:
\begin{eqnarray}
 \phi_{{\rm o}4}(\rho - 2 {\rm p}_{\rm o}\theta, \sigma + 2{\rm
q}_{\rm o} \theta + 2(\frac{2\pi}{N}i)) & = & \varsigma \left (
\frac {\rho - 2 {\rm p}_{\rm o}\theta }{2} \right ) \varsigma
\left ( \frac {\rho}{2} \right ) \varsigma \left ( \frac {\sigma +
2{\rm q}_{\rm o} \theta +
2(\frac{2\pi}{N}i)}{2}  \right ) \nonumber \\
&& \label{iklmjknbdhncjhdm} \varsigma \left ( \frac {\sigma}{2}
\right ) \phi_{{\rm o}4}(\rho, \sigma). \end{eqnarray}

\noindent

A fundamental region for $\Omega_{\phi_{{\rm o}4}}$ which is
spanned by the  $\;{\rm S}^{1}(-{\rm p}_{{\rm o}},{\rm q}_{{\rm
o}})-{\rm orbits}\;$ is the region $\;{\rm F}_{4}\;$ which is
depicted in Figure \ref {bolos}.
The fundamental regions $\;{\rm
F}_{4}\;$ and $\;{\rm F}_{3}\;$ are symmetric with respect to the
axis $\;\sigma=-\rho.\;$

We now make the following change of co$-$ordinates: \be
\label{iouryyyrueikjuhyhf} u_{4}=-\rho - {\rm q} \sigma, \qquad
v_{4}= \sigma-{\rm q} \rho, \quad \rm{where} \ {\rm q}\equiv \frac {{\rm
p}_{\rm o}}{{\rm q}_{\rm o}}. \ee The new coordinate $u_{4}$
remains constant along each \ ${\rm S}^{1}_{(- {\rm p}_{\rm o},
{\rm q}_{\rm o})}$$-$orbit. The lines of the two families of lines
$u_{4}$=constant and $v_{4}$=constant are straight and
perpendicular to each other. In the $ u_{4},v_{4} $ coordinates
the fundamental region we chose is given by \be 0 \leq u_{4} \leq
\frac {2 \pi}{{\rm q}_{\rm o}}, \qquad 0 \leq v_{4} \leq  {2
\pi}{\rm q}_{\rm o}(1+{\rm q}^{2}). \ee We define \be
\label{hdddjnbhgvbdgftry} \psi_{{\rm o}4}(u_{4},v_{4})=\phi_{{\rm
o}4} \left ( -\frac {u_{4}+{\rm q}v_{4}}{1+{\rm q}^{2}}, \frac
{v_{4}-{\rm q}u_{4}} {1+{\rm q}^{2}} \right ) \equiv \phi_{{\rm
o}4}(\rho, \sigma).  \ee The function $ \; \psi_{{\rm
o}4}(u_{4},v_{4}) \;$ satisfies the following periodicity
condition:   \be \label{kilokiiiijuhynbv} \psi_{{\rm o}4}(u_{4}+a
\frac{2\pi}{{\rm q}_{{\rm o}}} \ , \ v_{4}+ b \frac{2\pi}{{\rm
q}_{{\rm o}}}) = \psi_{{\rm o}4}(u_{4},v_{4}),  \ee where $\; a
\;$ and $\; b \;$ are integers. This last periodicity condition
is equivalent to the following periodicity conditions: \be
\label{jkinnbbhdmknjhbh} \psi_{{\rm o}4}(u_{4} - 2 \pi{\rm q},
v_{4}+2\pi)=\psi_{{\rm o}4}(u_{4}, v_{4}), \quad \quad \psi_{{\rm
o}4}(u_{4}-2\pi, v_{4} - 2 \pi{\rm q})=\psi_{{\rm o}4}(u_{4},
v_{4}).  \ee

\noindent The functional Eq (\ref {iklmjknbdhncjhdm}) in the new
co$-$ordinates $u_{4}, \ v_{4}$ is written as follows:
\begin{eqnarray}
\psi_{{\rm o}4}(u_{4}-{\rm q} \frac{4 \pi}{N}\nu ,v_{4}+(1-{\rm
q}^{2})2{\rm q}_{\rm o}\theta+ \frac{4\pi}{N}\nu) & = & \varsigma
\left ( -\frac {u_{4}+{\rm q}v_{4} + 2 {\rm p}_{\rm o} \theta (1 +
{\rm q}^{2}) } { 2 (1 + {\rm q}^{2})} \right ) \varsigma \left (
-\frac {u_{4}+{\rm q}v_{4}}{2(1+{\rm q}^{2})} \right )  \nonumber
\\ &  & \varsigma \left ( \frac {v_{4}- {\rm q} u_{4} + 2{\rm q}_
{\rm o}  \theta (1 + {\rm q}^{2})+ \frac{4 \pi}{N}\nu(1+{\rm
q}^{2}) }  { 2 (1 + {\rm q}^{2})} \right ) \nonumber \\ &&
\label{iiijjhggdjkfmnj} \varsigma \left ( \frac
 {v_{4} - {\rm q}u_{4}}
{2(1+{\rm q}^{2})} \right ) \psi_{{\rm o}4}(u_{4},v_{4}).
\end{eqnarray}
By setting  $\; v_{4}=0, \;$  $\; \nu =0 \;$ and by
 restricting $ \; u_{4} \;$ to the range of
values $ \; u_{4} \in [0,\frac{2\pi}{{\rm q}_{{\rm o}}}], \;$ we
obtain from Eq. (\ref{iiijjhggdjkfmnj}) \be
\label{gbbvcsbbnbanahhujmmsl} \psi_{{\rm o}4}(u_{4},v_{4})=
\varsigma \left ( \frac { u_{4} +{\rm q} v_{4} }{2(1+{\rm q}^{2})}
\right ) \varsigma \left ( \frac
 { {\rm q}u_{4} - v_{4} }
{2(1+{\rm q}^{2})} \right )  \psi_{{\rm o}4}(u_{4},0).  \ee The
function $ \; \psi_{{\rm o}4}(u_{4},0) \;$ satisfies  the
following periodicity condition:
\be\label{kbbcnnnncmjdkjuirrrrjhf} \psi_{{\rm
o}4}(u_{4}+\frac{2\pi}{{\rm q}_{{\rm o}}N}l,0)=-\psi_{{\rm
o}4}(u_{4},0),  \ee where $\;l={\rm gcd}(N,{\rm p}_{{\rm
o}}).\;$ Eq. (\ref{kbbcnnnncmjdkjuirrrrjhf}) is equivalent to:
\be\label{iikjumnhjkdkolmnjfhg} \psi_{{\rm
o}4}(u_{4}+\frac{2\pi}{{\rm q}_{{\rm o}}N'},0)=-\psi_{{\rm
o}4}(u_{4},0),  \ee where $\;N=N'l,\;$ $\;{\rm p}_{{\rm o}}={\rm
p}_{{\rm o}}'l,\;$ and $\;N', \ {\rm p}_{{\rm o}}'\;$ are coprime.
From Eqs. (\ref {iouryyyrueikjuhyhf}), (\ref {hdddjnbhgvbdgftry})
and (\ref {gbbvcsbbnbanahhujmmsl}) we conclude that the solution
to the functional Eq.(\ref {iklmjknbdhncjhdm}) is given almost
everywhere by: \be \label{hdnbmjknjikolmknjhds} \phi_{{\rm
o}4}(\rho,\sigma)= \varsigma \left ( \frac {\sigma}{2} \right )
\varsigma \left ( \frac {\rho}{2} \right ) \vartheta_{{\rm q}N'}
(-\rho - {\rm q} \sigma), \ee where $\vartheta_{{\rm q}N'} (u_{4})
\equiv \psi_{{\rm o}4}(u_{4},0)$.
\noindent The square$-$integrable condition on $\;\phi_{{\rm
o}4}(\rho,\sigma)\;$ gives: \be \int _{0}^{2\pi / {\rm q}_{\rm
o}N'} | \vartheta_{{\rm q}N'} (u_{4}) |^ {2} {\rm d}u_{4} < +
\infty, \ee Therefore, $A(H(N, -{\rm p}_{\rm o}, {\rm q}_{\rm
o}))$, as a Hilbert space, is isomorphic to ${\rm L}^{2}_{4{\rm
q}N'}({\rm A}_{4{\rm q}_{{\rm o}}N'}, \lambda, R)$, the Hilbert
space of the real valued functions $ \vartheta _{{\rm
q}N'}(u_{4})$ which are defined on the arc ${\rm A}_{4{\rm
q}_{{\rm o}}N'}$ $(v_{4}=0, \; 0 \leq u_{4} \leq \frac {2
\pi}{{\rm q}_{\rm o}N'})$ and are square$-$integrable with respect
to the Lebesgue measure $\; \lambda \;$ along this arc. This
completes the proof. \vspace*{4mm}
\end{pf}

\noindent
Note that in Appendix \ref{s6} we find the elementary regions of the actions 
of the groups ${\rm S}^{1}_{( {\rm q}_{\rm o},
{\rm p}_{\rm o})}$ 
on the torus $\rm{T} \simeq P_{1}(R) \times P_{1}(R)$
in the various cases. From these one can easily 
find the corresponding elementary regions of the actions
of the groups $H(N,{\rm q}_{\rm o},{\rm
p}_{\rm o})\simeq C_{N} \times {\rm S}^{1}_{( {\rm q}_{\rm o},
{\rm p}_{\rm o})}$ on the torus $\rm{T} \simeq P_{1}(R) \times P_{1}(R)$.

 \noindent
From Eq.(\ref{tzatziki}) we have $\;
(T'(-I,-I)\phi)(x,y)=\phi(x,y). \;$ Therefore, every little group
must contain the element $\;(-I,-I).\;$ Moreover, a group
$\;S\subseteq SO(2) \times SO(2)\;$ is a little group if and only
if its invariant space of functions $\;A(S)\;$ contains an element
$\; \phi_{{\rm o}}\;$ which is invariant under $\; S \;$ and is
not invariant under any bigger group $\; S' \supset S. \;$ All the
groups $\;H(N,{\rm q},{\rm p}),\;$where $\;N\;$ is odd and $\;{\rm
q},{\rm p}\;$ are odd and relatively prime, contain the element
$\;(-I,-I),\;$ (take $\;\nu=0\;$ and $\; \theta=\pi).\;$What
remains to be proved
is
that for every group $\;H(N,{\rm q},{\rm p})\;$ the associated
space of invariant vectors $\;A(H(N,{\rm q},{\rm p}))\;$ contains
at least one element $\; \phi_{{\rm o}}\;$ which is invariant
under $\;H(N,{\rm q},{\rm p})\;$ and is not invariant under any
bigger group $\;S',\;$ where, $\;H(N,{\rm q},{\rm p})\subset
S'\subseteq SO(2) \times SO(2).\;$ For this reason we call the
groups $\;H(N,{\rm q},{\rm p})\;$ potential little groups.
Summarizing the previous results we have the following theorem.



\noindent
\begin{thrm}
The infinite potential little groups of $\;\H \B
\;$ are the
one$-$dimensional groups $\;H(N,{\rm q},{\rm p}),\;$ where $\;N\;$
is odd and $\;{\rm q}, \ {\rm p}\;$ are odd and relatively prime.
\end{thrm}

\section{All Infinite Potential Little Groups are Actual}

\label{s3}

We proceed now to carry out step (III) of the Programme. Most
likely most functions in $\;A(H(N,{\rm q},{\rm p}))\;$ are
invariant only under $\;H(N,{\rm q},{\rm p})\;$ and do not have
any higher invariance. {\it A priori} \normalfont we can not
exclude the possibility that for some group(s) $\;H'(N,{\rm
q},{\rm p})\;$ {\it all} \normalfont the functions in
$\;A(H'(N,{\rm q},{\rm p}))\;$ are also invariant under some
bigger group(s) $\; S , \;$ where, $\;H'(N,{\rm q},{\rm p})
\subset S \subseteq SO(2) \times SO(2). \;$ If this were the case,
then these groups $\;H'(N,{\rm q},{\rm p})\;$ would not qualify as
(actual) little groups. If for some group $\;H(N,{\rm q},{\rm
p})\;$ there is only one element $\; \phi_{{\rm o}} \in A(H(N,{\rm
q},{\rm p}))\;$ which is invariant only under $\;H(N,{\rm q},{\rm
p})\;$ and which has no higher invariance, then the potential
little group $\;H(N,{\rm q},{\rm p})\;$ does qualify as an actual
little group. Therefore, in order to prove that a potential little
group $\;H(N,{\rm q},{\rm p})\;$ is an actual little group, it
suffices to prove that there exists an element $\; \phi_{{\rm o}}
\in A(H(N,{\rm q},{\rm p})) \;$ which has no higher invariance. It
will be shown that {\it every} potential little group $\;H(N,{\rm
q},{\rm p}),\;$ is an actual little group. The proof is based on
explicit construction. For {\it every} \normalfont $\;A(H(N,{\rm
q},{\rm p})),\;$ a function $\; \phi_{{\rm o}} \in A(H(N,{\rm
q},{\rm p})) \;$ is constructed which is invariant only under
$\;H(N,{\rm q},{\rm p})\;$ and which has no higher invariance. We
will restrict attention to the case 1. The treatment of cases 2, 3,
and 4 is similar (the cases 1, 2, 3, and 4 are defined in Theorem
\ref{iiiuolokjuuuihyjuik}). The details are as follows.

\noindent For a given $\;H(N,{\rm q}_{{\rm o}},{\rm p}_{{\rm
o}})\;$
we define the following element of $\;A(H(N,{\rm
q}_{{\rm o}},{\rm p}_{{\rm o}})) \;$ \be \label{eerrregfdhjunbg}
\psi_{{\rm o}1}(u_{1},v_{1})= \varsigma \left ( \frac { {\rm
q}u_{1} - v_{1} }{2(1+{\rm q}^{2})} \right ) \varsigma \left (
\frac
 { u_{1} + {\rm q}v_{1} }
{2(1+{\rm q}^{2})} \right )  \sin \left( \frac{N{\rm q}_{{\rm
o}}}{2l}u_{1}\right) ,  \ee where the function $\;\sin \left(
\frac{N{\rm q}_{{\rm o}}}{2l}u_{1}\right)\equiv\psi_{{\rm
o}1}(u_{1},0)\;$ satisfies the periodicity condition $\;
\psi_{{\rm o}1}(u_{1}+\frac{2\pi}{{\rm q}_{{\rm
o}}N}l,0)=-\psi_{{\rm o}1}(u_{1},0),
\;$(Eq.(\ref{mikiloooopoijijik})). It will be shown that if $k\in
{\rm S}{\rm O}(2) \times {\rm S}{\rm O}(2)$ and $\;\left (
T^{\prime }(k)\psi_{{\rm o}1} \right )(u_{1},v_{1})=\psi_{{\rm
o}1}(u_{1},v_{1})\;$, then $\;k\in H(N,{\rm q}_{{\rm o}},{\rm
p}_{{\rm o}})\;$. \ Here, equality is in the Hilbert space sense.
\ In other words, it will be shown that
\begin{equation}
\label{iiiiuuunnbhdseruikolko} \left\| \left ( T^{\prime
}(k)\psi_{{\rm o}1} \right )(u_{1},v_{1})-\psi_{{\rm
o}1}(u_{1},v_{1})\right\| =0\quad \Longrightarrow \quad k\in
H(N,{\rm q}_{{\rm o}},{\rm p}_{{\rm o}}),
\end{equation}
where $\;\psi_{{\rm o}1}(u_{1},v_{1})\;$ is given by Eq.
(\ref{eerrregfdhjunbg}). Let $\; k = \left(  \left(
\begin{array}{cc}
\cos\omega& \sin\omega \\
\!\!\!\!-\sin\omega & \cos\omega
\end{array}
\right),
 \left( \begin{array}{cc}
\cos\chi & \sin\chi \\
\!\!\!\!\!-\sin\chi & \cos\chi
\end{array}
\right) \right)\equiv \left( R(\omega),R(\chi) \right). \;$ For
{\it any} \normalfont $\; \psi_{{\rm o}1}(u_{1},v_{1})\;$ the
expression $\;\left ( T^{\prime }(k)\psi_{{\rm o}1} \right
)(u_{1},v_{1})\;$ reads
\begin{eqnarray} \left ( T^{\prime }(k)\psi_{{\rm o}1} \right
)(u_{1},v_{1}) &=& \varsigma \left ( \frac{v_{1}-{\rm
q}u_{1}+2\omega(1+{\rm q}^{2})}{2(1+{\rm q}^{2})}\right)\varsigma
\left ( \frac{v_{1}-{\rm
q}u_{1}}{2(1+{\rm q}^{2})}\right) \nonumber \\
&&\varsigma \left ( \frac{v_{1}{\rm q}+u_{1}+2\chi(1+{\rm
q}^{2})}{2(1+{\rm q}^{2})}\right)\varsigma \left ( \frac{v_{1}{\rm
q}+u_{1}}{2(1+{\rm q}^{2})}\right) \nonumber \\
\label{iiioukjhfuydiol} &&\psi_{{\rm o}1}(u_{1}+2(\chi-{\rm
q}\omega),v_{1}+2(\omega+{\rm q}\chi)).
\end{eqnarray}
\noindent
For $\; \psi_{{\rm o}1}(u_{1},v_{1})= \varsigma \left ( \frac {
{\rm q}u_{1} - v_{1} }{2(1+{\rm q}^{2})} \right ) \varsigma \left
( \frac
 { u_{1} + {\rm q}v_{1} }
{2(1+{\rm q}^{2})} \right )   \sin \left( \frac{N{\rm q}_{{\rm
o}}}{2l}u_{1}\right) \;$ the expression $\;\left ( T^{\prime
}(k)\psi_{{\rm o}1} \right )(u_{1},v_{1})\;$
(Eq.(\ref{iiioukjhfuydiol})) reads
\begin{eqnarray} \left ( T^{\prime }(k)\psi_{{\rm o}1} \right
)(u_{1},v_{1}) &=& \varsigma \left ( \frac{{\rm
q}u_{1}-v_{1}}{2(1+{\rm q}^{2})}\right) \varsigma \left (
\frac{v_{1}{\rm
q}+u_{1}}{2(1+{\rm q}^{2})}\right) \nonumber \\
\label{ujkhhhnbghytvvbdfgfyujnm} &&\sin \left [ \frac{N{\rm
q}_{{\rm o}}}{2l}(u_{1}+2(\chi-{\rm q}\omega)) \right ].
\end{eqnarray}
Taking into account expressions (\ref{eerrregfdhjunbg}) and
(\ref{ujkhhhnbghytvvbdfgfyujnm}) the condition (Eq.
(\ref{iiiiuuunnbhdseruikolko})) $\; \left\| \left ( T^{\prime
}(k)\psi_{{\rm o}1} \right ) \right . \;$\newline $\; \left .
(u_{1},v_{1})-\psi_{{\rm o}1}(u_{1},v_{1})\right\| =0 \;$ reads
\begin{eqnarray}
\int_{{\rm F}_{1}} \int \left [ \left ( T^{\prime }(k)\psi_{{\rm
o}1} \right ) (u_{1},v_{1})-\psi_{{\rm o}1}(u_{1},v_{1}) \right
]^{2} |J| {\rm d}u_{1} \, {\rm d}v_{1} & = & 0 \; \Longrightarrow
\nonumber \\
\int_{{\rm F}_{1}} \int  \left ( \sin \left [ \frac{N{\rm q}_{{\rm
o}}}{2l}(u_{1}+2(\chi-{\rm q}\omega)) \right ]- \sin \left [
\frac{N{\rm q}_{{\rm o}}}{2l}u_{1}\right ] \right )^{2}
 |J| {\rm d}u_{1} \, {\rm d}v_{1} & = & 0 \; \Longrightarrow
\nonumber \\
\label{uiiijhyughdgryuio} \int_{{\rm F}_{1}} \int  \left [ \sin
\left ( \frac{N{\rm q}_{{\rm o}}}{2l}u_{1} \right ) \left \{ \cos
\left ( \frac{N{\rm q}_{{\rm o}}}{2l}2(\chi-{\rm q}\omega) \right
) -1 \right \} \right . & + & \nonumber \\ \left . \cos \left (
\frac{N{\rm q}_{{\rm o}}}{2l}u_{1} \right ) \sin \left (
\frac{N{\rm q}_{{\rm o}}}{2l}2(\chi-{\rm q}\omega) \right ) \right
]^{2}
 |J| {\rm d}u_{1} \, {\rm d}v_{1} & = & 0,
\end{eqnarray}
where ${\rm F}_{1}$ denotes the fundamental region depicted in
Fig. (\ref{melas}) and $\;|J|\;$ denotes the absolute value of
$\;J,\;$ the Jacobian $ \left | \frac { \partial (\rho, \sigma)}
{\partial (u_{1}, v_{1})} \right |  $ of the co-ordinate
transformation $ (u_{1},v_{1}) \longmapsto (\rho, \sigma)$. After
some simple algebra Eq. (\ref{uiiijhyughdgryuio}) gives \be
\label{uuuiiiiokilojhyhfgrt}2\pi \left [ \pi \left \{ \cos \left (
\frac{N{\rm q}_{{\rm o}}}{l}(\chi-{\rm q}\omega) \right ) -1
\right \}^{2} + \pi \left \{\sin \left ( \frac{N{\rm q}_{{\rm
o}}}{l}(\chi-{\rm q}\omega) \right ) \right \}^{2}\right ]=0 .\ee

\noindent The sum of squares in (\ref{uuuiiiiokilojhyhfgrt}) can
only vanish if each individual square vanishes, so
(\ref{uuuiiiiokilojhyhfgrt}) vanishes if and only if $\;\cos \left
( \frac{N{\rm q}_{{\rm o}}}{l}(\chi-{\rm q}\omega) \right ) =1\;$
and $\;\sin \left ( \frac{N{\rm q}_{{\rm o}}}{l}(\chi-{\rm
q}\omega) \right )=0 .\;$ \ Therefore, for some integer $\;a,\;$
we must have
\begin{equation}
\label{iwerdfgtyhsjau} \frac{N{\rm q}_{{\rm o}}}{l}(\chi-{\rm
q}\omega) =2\pi a \quad {\rm or} \quad 2(\chi-{\rm
q}\omega)=\frac{l}{N{\rm q}_{{\rm o}}}4\pi a .
\end{equation}
Let $\;l\;$ be the highest common factor of $\;N\;$ and $\;{\rm
q}_{{\rm o}},\;$ and let, $\;N'\;$ and $\;{\rm q}_{{\rm o}}'\;$ be
the coprime numbers which are defined by $\;N=N'l\;$ and $\;{\rm
q}_{{\rm o}}={\rm q}_{{\rm o}}'l.\;$ Let $\;b\;$ and $\;\nu\;$ be
integer numbers. Since $\;N'\;$ and $\;{\rm q}_{{\rm o}}'\;$ are
coprime the difference $\;{\rm q}_{{\rm o}}'\nu-bN'=a\;$ takes all
possible integer values $\;a\;$ as the integers $\;b\;$ and
$\;\nu\;$ vary. The equality $\;a+bN'={\rm q}_{{\rm o}}'\nu\;$ is
equivalent to \be \label{iiiolplkiojuhjuy}4\pi a \frac{l}{N{\rm
q}_{{\rm o}}} + 2b\frac{2\pi}{{\rm q}_{{\rm o}}}=\frac{4\pi}{N}\nu
.\ee Therefore, for {\it every} \normalfont $\;a\;$ there always
exist integers $\;b\;$ and $\;\nu\;$ such that Eq.
(\ref{iiiolplkiojuhjuy}) is satisfied. Now according to Eq.
(\ref{iiioukjhfuydiol}) under the action of the element $\;k
\equiv (R(\omega),R(\chi)) \in {\rm S}{\rm O}(2) \times {\rm
S}{\rm O}(2)\;$ the points $\;(u_{1},v_{1})\;$ of the 2$-$torus
$\;{\rm T}^{2}\simeq P_{1}(R) \times P_{1}(R) \;$ are sent to \be
\label{iiiodfgsretbvfdg} u_{1} \longmapsto
u_{1}'=u_{1}+2(\chi-{\rm q}\omega) \quad , \quad v_{1}\longmapsto
v_{1}'=v_{1}+2(\omega+{\rm q}\chi). \ee The periodicity
condition (\ref{koiiuiuyhyjuik}) shows that the coordinate
$\;u_{1}\;$ is defined modulo $\; \frac{2\pi}{{\rm q}_{{\rm o}}}.
\;$ From Eqs. (\ref{iiiodfgsretbvfdg}) and (\ref{iwerdfgtyhsjau})
we obtain \be \label{aaacvfdgbhyru}u_{1}'-u_{1}=2(\chi-{\rm
q}\omega)=\frac{l}{N{\rm q}_{{\rm o}}}4\pi a .\ee By taking into
account Eq. (\ref{iiiolplkiojuhjuy}) we conclude that
(\ref{aaacvfdgbhyru}) is equivalent, mod. $\; \frac{2\pi}{{\rm
q}_{{\rm o}}},\;$ to \be \label{uuytrgdfrety} 2(\chi-{\rm
q}\omega)=\frac{4\pi}{N}\nu.  \ee By substituting \be
\label{ddsfgebhnjhy} \omega={\rm q}_{{\rm o}}\theta,   \ee
where $\; \theta \in R ,\;$ into Eq. (\ref{uuytrgdfrety}) we
obtain \be \chi={\rm p}_{{\rm o}}\theta+\frac{2\pi}{N}\nu.
\ee But these values of $\;\omega \;$ and $\;\chi \;$ are
precisely those for which $\;k\in H(N,{\rm q}_{{\rm o}},{\rm
p}_{{\rm o}}),$ so the $\;\psi _{{\rm o}1}(u_{1},v_{1})\;$ given
by (\ref{eerrregfdhjunbg}) does indeed have little group $\;L(\psi
_{{\rm o}1})=H(N,{\rm q}_{{\rm o}},{\rm p}_{{\rm o}}).$

\vspace{0.3cm} \noindent The previous results are summarized as
follows:
\begin{thrm}
All the infinite potential little groups of $\;\H \B
\;$ are
actual. These are the one$-$dimensional groups $\;H(N,{\rm q},{\rm
p}),\;$ where $\;N\;$ is odd and $\;{\rm q}, \ {\rm p}\;$ are odd and
relatively prime.
\end{thrm}

\section{Form of the induced representations}
\label{s4}

 \indent


Let $\; A \;$ and $\; \mathcal G \;$ be topological groups, and
let $\; T \;$ be a given homomorphism from $\; \mathcal G \;$ into
the group of automorphisms $\; {\rm Aut} (A) \;$ of $\; A. $
Suppose $\; A \;$ is abelian and $\; \mathcal H= A
\bigcirc\!\!\;\!\!\!\;\!\!\!\!s \  _{T} \mathcal G \;$ is the
semi$-$direct product of $\; A \;$ and $\; \mathcal G, $ specified
by the continuous action $\; T \; : \mathcal G \longrightarrow
{\rm Aut} (A). \;$ In the product topology of $\; A \times
\mathcal G \;$, $\; \mathcal H \;$ then becomes a topological
group. It is assumed that it becomes a separable locally compact
topological group.

\noindent
 In order to give the operators of the induced
representations explicitly it is necessary
(\cite{Wigner}, \cite{Mackey}, \cite{Mackey1}, \cite{Simms},
 \cite{Isham} and references therein) to give the following
information
\begin{enumerate}
\item{ An irreducible unitary representation $\;U\;$ of
$\;L_{\phi_{\rm{o}}}\;$ on a Hilbert space $\;D\;$ for each
$\;L_{\phi_{\rm{o}}}.\;$} \item{ A $\; \mathcal G\;$$-$quasi$-$invariant
measure $\; \mu \;$ on each orbit $\; \mathcal G \phi \approx
\mathcal G / L_{\phi_{\rm{o}}}; \;$ where $\; L_{\phi_{\rm{o}}} \;$ denotes the
little group of the base point $\; \phi_{\rm{o}} \in
A^{'} \;$ of the
orbit $\; \mathcal G \phi_{\rm{o}}; $ $A^{'}$ is the topological dual of $A$. }
\end{enumerate}

\noindent Let $\; D_{\mu} \;$ be the space of functions $\; \psi:
\mathcal G \rightarrow D \;$ which satisfy the conditions

\begin{eqnarray}
(a) & \psi(gl)=U(l^{-1})\psi(g) & (g \in \mathcal G, \ l \in
L_{\phi}) \nonumber \\
& & \nonumber \\ (b) & \int_{\mathcal {G} \phi_{\rm{o}}}
<\psi(q),\psi(q)>{\rm d}\mu(q) < \infty,  \nonumber
\end{eqnarray}
where the scalar product under the integral sign is that of $\; D.
\;$ Note, that the constraint $\; (a) \;$ implies that $\;
<\psi(gl),\psi(gl)>=<\psi(g),\psi(g)>, \;$ and therefore the inner
product $\;<\psi(g),\psi(g)>,$ $\ g \in \mathcal G, \;$ is constant
along  every element $\; q \;$ of the coset space $\; \mathcal
G/L_{\phi_{\rm{o}}} \approx \mathcal G \phi_{\rm{o}}.$ This allows to assign a
meaning to $\; <\psi(q),\psi(q)>, $  where $\; q=gL_{\phi_{\rm{o}}}, \;$ by
defining $\; <\psi(q),\psi(q)>:= <\psi(g),\psi(g)>. \;$ Thus the
integrand in $\; (b) \;$ becomes meaningful due to the condition
$\; (a). \;$ A pre$-$Hilbert space structure can now be given to $\;
D_{\mu} \;$ by defining the scalar product
\begin{equation}
\label{hujkiolokijuhy} <\psi_{1},\psi_{2}>=\int_{\mathcal G
\phi_{\rm{o}}}<\psi_{1}(q),\psi_{2}(q)> {\rm d}\mu(q),
\end{equation}
where $\; \psi_{1},\psi_{2} \in D_{\mu}. \;$ It is convenient to
complete the space $\; D_{\mu} \;$ with respect to the norm
defined by the scalar product (\ref{hujkiolokijuhy}). In the
resulting Hilbert space, functions are identified whenever they
differ, at most, on a set of $\; \mu-\!$measure zero. Thus our
Hilbert space is
\begin{equation}
 D_{\mu}=L^{2}(\mathcal G \phi_{\rm{o}}, \mu ,D).
\end{equation}
Define now an action of $\; \mathcal H= A
\bigcirc\!\!\;\!\!\!\;\!\!\!\!s \  _{T} \mathcal G \;$ on $\;
D_{\mu}\;$ by
\begin{eqnarray}
\label{e1}
(g_{{\rm o}}\psi)(q) & = & \sqrt{\frac{\rm{d}\mu_{g_{{\rm o}}}}{\rm{d}\mu}(q)} \psi(g_{{\rm o}}^{-1}q), \\
\label{e2}
\alpha \psi(q) &=&  e^{i <g_{{\rm o}} \phi_{\rm{o}},\alpha> } \psi(q)
\end{eqnarray}
where, $\; g_{{\rm o}} \in \mathcal G, \;$ $q \in \mathcal G \phi_{\rm{o}},$ and 
$\alpha \in A$. 
Eqs. (\ref{e1}) and (\ref{e2})   define the IRs of $\H \B$ induced for 
{\it each} $\; \phi_{\rm{o}} \in
A^{'} \;$
and {\it each} irreducible representation $\; U \;$ of
$\;L_{\phi_{\rm{o}}}.\;$ 
The `Jacobian' 
$\frac{\rm{d}\mu_{g_{{\rm o}}}}{\rm{d}\mu} $ 
of the group transformation is known as the \it Radon$-$Nikodym \normalfont derivative of 
$\mu_{g_{{\rm o}}}$ with respect to $\mu$ and ensures that the resulting IRs of $\H \B$ are unitary.


\noindent
 The central results of induced representation theory
(\cite{Wigner}, \cite{Mackey}, \cite{Mackey1}, \cite{Simms},
\cite{Isham} and references therein) are the following
\begin{enumerate}
\item {Given the topological restrictions on $\; \mathcal H= A
\bigcirc\!\!\;\!\!\!\;\!\!\!\!s \  _{T} \mathcal G \;$
(separability and local compactness), any representation of $\;
\mathcal H, \;$  constructed by the method above,  is irreducible
if the representation $\; U \;$ of $\; L_{\phi_{\rm{o}}} \;$ on $\; D \;$
is irreducible. Thus an irreducible representation of $\;
\mathcal H \;$ is obtained for {\it each} $\; \phi_{\rm{o}} \in
A^{'} \;$
and {\it each} irreducible representation $\; U \;$ of
$\;L_{\phi_{\rm{o}}}.\;$} \item {If $\; \mathcal H= A
\bigcirc\!\!\;\!\!\!\;\!\!\!\!s \  _{T} \mathcal G \;$ is a
regular semi$-$direct product (i.e., $\;
A^{'} \;$ contains a Borel
subset which meets each orbit in $\;
A^{'} \;$ under $\; \mathcal
H \;$ in just one point) then {\it all} of its irreducible
representations can be obtained in this way}.
\end{enumerate}

To find the induced representations of $\; \mathcal H \mathcal
B
= 
L^{2}(\mathcal P, \lambda, R)
\bigcirc\!\!\;\!\!\!\;\!\!\!\!\!s \  _{T}
\mathcal G, \;$  then,  it is enough to provide the information
cited in 1 and 2 for each of the orbit types.

\vspace{0.5cm} 1. \hspace{0.5cm} It was shown that little groups
are the groups $\; H(N,{\rm q}_{{\rm o}},{\rm p}_{{\rm o}})=C_{N}
\times {\rm S}^{1} _{({\rm q}_{{\rm o}},{\rm p}_{{\rm o}})}, \;$
where the numbers $\;N,{\rm q}_{{\rm o}},{\rm p}_{{\rm o}} \;$ are
odd, and moreover, the numbers $\; {\rm q}_{{\rm o}},{\rm p}_{{\rm
o}}\;$ are relatively prime. The little groups are abelian. All
IRs
of an abelian group are
one$-$dimensional. Firstly, we comment on the 
IRs
of $\;{\rm
S}^{1} _{({\rm q}_{{\rm o}},{\rm p}_{{\rm o}})}. \;$ Let $\;
R_{({\rm q}_{{\rm o}},{\rm p}_{{\rm o}})} \;$ be a complex 
one$-$dimensional representation of $\; {\rm S}^{1} _{({\rm q}_{{\rm
o}},{\rm p}_{{\rm o}})}. \;$ 
Then we may write \be R_{({\rm
q}_{{\rm o}},{\rm p}_{{\rm o}})}(\theta)= \chi_{({\rm q}_{{\rm
o}},{\rm p}_{{\rm o}})} (\theta) I, \ee where, $\; \theta
\rightarrow \chi_{({\rm q}_{{\rm o}},{\rm p}_{{\rm o}})} (\theta)
\;$ is a complex$-$valued function on $\; {\rm S}^{1} _{({\rm
q}_{{\rm o}},{\rm p}_{{\rm o}})} \;$ that is never zero, and $\;
I \;$ is the identity operator in a one$-$dimensional complex
Hilbert space $\; D \approx C \;$ ($\; C \;$=complex numbers).
Since $\; R_{({\rm q}_{{\rm o}},{\rm p}_{{\rm o}})} \;$ is a
representation \be \label{luhgvhm,gbvk,yutfvikt} \chi_{({\rm
q}_{{\rm o}},{\rm p}_{{\rm o}})} (\theta_{1} + \theta_{2} )=
\chi_{({\rm q}_{{\rm o}},{\rm p}_{{\rm o}})} (\theta_{1})
\chi_{({\rm q}_{{\rm o}},{\rm p}_{{\rm o}})} (\theta_{2}). \ee The
condition for the representations being unitary reads \be
\label{obvuklhbvuklygfvlybvk} | \chi_{({\rm q}_{{\rm o}},{\rm
p}_{{\rm o}})} (\theta) |=1. \ee The condition $\; {\rm S}^{1}
_{({\rm q}_{{\rm o}},{\rm p}_{{\rm o}})} (\theta + 2 \pi)= {\rm
S}^{1} _{({\rm q}_{{\rm o}},{\rm p}_{{\rm o}})} (\theta) \;$
implies \be \label{lyihbj,hgvyugfukygvlukygv} \chi_{({\rm q}_{{\rm
o}},{\rm p}_{{\rm o}})} (\theta+2 \pi)= \chi_{({\rm q}_{{\rm
o}},{\rm p}_{{\rm o}})} ( \theta ). \ee Therefore, one has to find
continuous complex$-$valued functions $\; \chi (\theta) \;$
satisfying the equations
(\ref{luhgvhm,gbvk,yutfvikt}), (\ref{obvuklhbvuklygfvlybvk}) and
(\ref{lyihbj,hgvyugfukygvlukygv}). It is well known (see for
example \cite{Vilenkin}, page 70) that all such functions have the
form

\be \chi_{({\rm q}_{{\rm o}},{\rm p}_{{\rm o}})}(\theta)= e^{i n
\theta}, \quad {\rm where} \; {\rm n} \; {\rm is} \; {\rm an} \;
{\rm integer}. \ee It is worth pointing out that $ \; \chi_{({\rm
q}_{{\rm o}},{\rm p}_{{\rm o}})}(\theta) \;$ does not depend on
the pair $ \; ({\rm q}_{{\rm o}},{\rm p}_{{\rm o}}).\;$
Consequently, the IRs of
 $
{\rm S}^{1}_{( {\rm q}_{\rm o}, {\rm p}_{\rm o})} \equiv
  \left(  \left [ \begin{array}{cc}
\cos({\rm q}_{\rm o}\theta) & \sin({\rm q}_{\rm o}\theta) \\
\!\!\!\!\!-\sin({\rm q}_{\rm o}\theta) & \cos({\rm q}_{\rm
o}\theta)
\end{array}
\right ],
 \left [ \begin{array}{cc}
\cos({\rm p}_{\rm o}\theta) & \sin({\rm p}_{\rm o}\theta) \\
\!\!\!\!\!-\sin({\rm p}_{\rm o}\theta) & \cos({\rm p}_{\rm
o}\theta)
\end{array}
\right ] \right) $ are indexed by an integer $\; n \;$ which for
distinct representations takes the values
$\;n=...,-2,-1,0,1,2,...\;$ and are given by multiplication in one
complex dimension $\;D_{1} \approx C\;$ by
 \be
 D^{(n)}
\left (
  \left(  \left [ \begin{array}{cc}
\cos({\rm q}_{\rm o}\theta) & \sin({\rm q}_{\rm o}\theta) \\
\!\!\!\!\!-\sin({\rm q}_{\rm o}\theta) & \cos({\rm q}_{\rm
o}\theta)
\end{array}
\right ],
 \left [ \begin{array}{cc}
 \cos({\rm p}_{\rm o}\theta) & \sin({\rm p}_{\rm o}\theta) \\
\!\!\!\!\!-\sin( {\rm p}_{\rm o}\theta) & \cos({\rm p}_{\rm
o}\theta)
\end{array}
\right ] \right) \right )= e^{in\theta} . \ee We comment now on
the IRs of the cyclic group $\; C_{N} . \;$

\vspace{0.5cm} \noindent The IRs $\;U_{N} \;$ of the cyclic group
$\; C_{N} \;$ are well known (see for instance \cite{Mackey2}).
They are indexed by an integer $\; \nu \;$ which , for distinct
representations,takes values in the set $\; \nu \in
\{0,1,2,...,N-1 \} \;$. Thus, the number of 
IRs
of $\; C_{N} \;$ equals to the order
of the group. Denoting them by $\; D^{(\nu)} \;$, they are given
by multiplication in one complex dimension $\; D_{2} \approx C \;$
by
\begin{equation}
\label{;ujujuunlwjq;jnqj} D^{(\nu)}\left( \left ( {\rm I} , \left[
\begin{array}{cc}
\cos \frac{2\pi}{N}{\rm j} & \sin \frac{2\pi}{N}{\rm j} \\
-\sin \frac{2\pi}{N}{\rm j} & \cos \frac{2\pi}{N}{\rm j}
\end{array}
\right]
 \right )
 \right )=e^{i\frac{2 \pi}{N}\nu{\rm j}} ,
\end{equation}
where j parameterizes the elements of the group $\;C_{N}.\;$

\vspace{0.5cm} \noindent It follows
that the IRs of $\; H(N,{\rm q}_{{\rm o}},{\rm p}_{{\rm
o}})=C_{N} \times {\rm S}^{1}_{({\rm q}_{{\rm o}}, {\rm p}_{{\rm
o}})} \;$ are indexed by the pair of integers $\; (\nu,n). \;$ The
indices $\; \nu \;$ and $\; n \;$ for distinct 
IRs
of
 $\; H(N,{\rm q}_{{\rm o}},{\rm p}_{{\rm o}})\;$
 take independently values in the sets
$\;\{0,1,2,...,N-1\}\;$ and $\;{\rm Z},\;$ where $\;{\rm Z}\;$
denotes the set of integers. Denoting them by $\;D^{(\nu,n)} \;$,
they are given by multiplication in one complex dimension
$\;D\approx C \;$ by
\begin{eqnarray}
 D^{(\nu,n)}
\left (
  \left(  \left [ \begin{array}{cc}
\cos({\rm q}_{{\rm o}}\theta) & \sin({\rm q}_{{\rm o}}\theta) \\
\!\!\!\!\!-\sin({\rm q}_{{\rm o}}\theta) & \cos({\rm q}_{{\rm
o}}\theta)
\end{array}
\right ],
 \left [ \begin{array}{cc}
 \cos({\rm p}_{{\rm o}}\theta+\frac{2\pi}{N}{\rm j}) & \sin({\rm p}_{{\rm o}}\theta+\frac{2\pi}{N}{\rm j}) \\
\!\!\!\!\!-\sin({\rm p}_{{\rm o}}\theta+\frac{2\pi}{N}{\rm j}) &
\cos({\rm p}_{{\rm o}}\theta+\frac{2\pi}{N}{\rm j})
\end{array}
\right ] \right) \right ) & = &  \nonumber \\
\!\!\!\!\!\!\!\!\!\!\!\!\!e^{i\frac{2 \pi}{N}\nu{\rm j}}
e^{in\theta}. &&
\end{eqnarray}

\vspace{0.5cm} \noindent
 We now proceed to give the information cited in 2. Although a
 $\; \mathcal G$$-$quasi$-$invariant
measure  is all what is needed,  a $\; \mathcal G$$-$invariant
measure  will be provided in all cases.

\vspace{0.5cm} 2. \hspace{0.5cm} We want to  construct a $\;
\mathcal G$$-$invariant measure on the orbits $\; 01 \equiv \mathcal
G/ L_{\phi}= ({\rm S}{\rm L}(2,R) \times {\rm S}{\rm L}(2,R))/
H(N,{\rm q}_{{\rm o}},{\rm p}_{{\rm o}})=({\rm S}{\rm L}(2,R)
\times {\rm S}{\rm L}(2,R))/ ({\rm S}^{1}({\rm q}_{{\rm o}},{\rm
p}_{{\rm o}}) \times C_{N}) = (({\rm S}{\rm L}(2,R) \times {\rm
S}{\rm L}(2,R))/{\rm S}^{1}({\rm q}_{{\rm o}},{\rm p}_{{\rm
o}}))/C_{N}.\;$ It turns out that it suffices to construct
invariant measure on the orbits $\; \widetilde {01} \approx (({\rm
S}{\rm L}(2,R) \times {\rm S}{\rm L}(2,R))/{\rm S}^{1}({\rm
q}_{{\rm o}},{\rm p}_{{\rm o}}) .\;$

Firstly, we construct invariant measure on the orbits $\;
\widetilde {01} \equiv (({\rm S}{\rm L}(2,R) \times {\rm S}{\rm
L}(2,R))/{\rm S}^{1}({\rm q}_{{\rm o}},$ \newline $\; {\rm
p}_{{\rm o}}) .\;$ It will be shown in the Appendix
\ref{s7} 
that a $\; \mathcal G$$-$invariant measure
on the orbits $\; \widetilde {01} \equiv (({\rm S}{\rm L}(2,R)
\times {\rm S}{\rm L}(2,R))/{\rm S}^{1}({\rm q}_{{\rm o}},{\rm
p}_{{\rm o}}) \;$ may be constructed from an $\; \mathcal
G$$-$invariant measure on $\; \mathcal G \;$ and an $\;{\rm
S}^{1}({\rm q}_{{\rm o}},{\rm p}_{{\rm o}})-$invariant measure on
$\;{\rm S}^{1}({\rm q}_{{\rm o}},{\rm p}_{{\rm o}}).$


We start by constructing an $\; \mathcal G$$-$invariant measure on
$\;\mathcal G=G \times G, \;$ $\; G={\rm S}{\rm L}(2,R). \;$
 A $\;G$$-$invariant measure  $\; \mu_{G} \;$ on $\; G={\rm S}{\rm
L}(2,R) \;$ is well known. Indeed, in \cite{gel}, p.214$-$215 an
invariant measure on $\;{\rm S}{\rm L}(2,C)\;$ is explicitly
constructed. With a similar construction one can obtain an $\; G$$-$
invariant measure on $\;{\rm S}{\rm L}(2,R).\;$ So, if ${\rm
S}{\rm L}(2,R)= \left \{ \left(
\begin{array}{cc}
a  & b \\
c   &  d
\end{array}
\right), \qquad a,b,c,d \in R, \quad ad-bc=1 \right \}, \;$ an
invariant measure on $\;{\rm S}{\rm L}(2,R)\;$ is given by \be
\label{lyhgbhgblyuitg7lyugk} {\rm d}g= \frac{{\rm d}a \wedge {\rm
d}b \wedge {\rm d}c}{a}. \ee $\; G$$-$invariant measure means
$\; {\rm d}g={\rm d}(gg_{{\rm o}})={\rm d} (g_{{\rm o}}g), \quad
g_{{\rm o}} \in {\rm S}{\rm L}(2,R). \;$ 
The group $\; \mathcal
G={\rm S}{\rm L}(2,R) \times {\rm S}{\rm L}(2,R) \;$ is the set of
matrices
$$
\mathcal G = \left \{  \left(  \left( \begin{array}{cc}
a  & b \\
c   &  d
\end{array}
\right),
 \left( \begin{array}{cc}
 e  & f  \\
 j  & k
\end{array}
\right) \right), \qquad a,b,c,d,e,f,j,k \in R, \quad ad-bc=1 ,
\quad ek-jf=1 \right \}.
$$
Using the fact that the measure $\;{\rm d}g, \;$
(Eq.\ref{lyhgbhgblyuitg7lyugk}), is invariant on $\;{\rm S}{\rm
L}(2,R),\;$ one can easily prove that the
measure $\; {\rm d}{\rm g} \;$ on $\; \mathcal G= {\rm S}{\rm
L}(2,R) \times {\rm S}{\rm L}(2,R)\;$
given by \be \label{kugvukgvklutyfv} {\rm d}{\rm g}=\frac{{\rm d}a
\wedge {\rm d}b \wedge {\rm d}c}{a} \bigwedge \frac{{\rm d}e
\wedge {\rm d}f \wedge {\rm d}j}{e} \ee is $\; \mathcal
G$$-$invariant measure on $\; \mathcal G. \;$ Thus, an invariant
measure on  $\; \mathcal G \;$ is obtained.

We proceed now to give an invariant measure on $\;S^{1}({\rm
q}_{{\rm o}},{\rm p}_{{\rm o}}).\;$ An $\; S^{1}({\rm q}_{{\rm
o}},{\rm p}_{{\rm o}})-$invariant measure on $\;S^{1}({\rm
q}_{{\rm o}},{\rm p}_{{\rm o}})\;$ is provided by the usual
Lebesgue measure $\; {\rm d} \theta \;$ on $\;S^{1}({\rm q}_{{\rm
o}},{\rm p}_{{\rm o}}).\;$

We obtained invariant measures on the groups $\; \mathcal G \;$
and $\; S^{1}({\rm q}_{{\rm o}},{\rm p}_{{\rm o}}).\;$ By using
these invariant measures one can construct a $\; \mathcal G$$-$ invariant 
measure on the orbits $\; \widetilde {01} \equiv (({\rm
S}{\rm L}(2,R) \times {\rm S}{\rm L}(2,R))/{\rm S}^{1}({\rm
q}_{{\rm o}}, {\rm p}_{{\rm o}}) .\;$ This construction is given
explicitly in the Appendix \ref{s7}.
The $\; \mathcal G$$-$invariant 
measure on the orbits $\; \widetilde {01} \equiv
(({\rm S}{\rm L}(2,R) \times {\rm S}{\rm L}(2,R))/{\rm S}^{1}({\rm
q}_{{\rm o}}, {\rm p}_{{\rm o}}) \;$ so constructed is also $\;
\mathcal G - $ invariant measure on the orbits $\; 01 \equiv
\mathcal G/ L_{\phi}= ({\rm S}{\rm L}(2,R) \times {\rm S}{\rm
L}(2,R))/ H(N,{\rm q}_{{\rm o}},{\rm p}_{{\rm o}})=({\rm S}{\rm
L}(2,R) \times {\rm S}{\rm L}(2,R))/ ({\rm S}^{1}({\rm q}_{{\rm
o}},{\rm p}_{{\rm o}}) \times C_{N}) = (({\rm S}{\rm L}(2,R)
\times {\rm S}{\rm L}(2,R)) / $ \newline ${\rm S}^{1}({\rm
q}_{{\rm o}},{\rm p}_{{\rm o}}))/C_{N}.\;$ We now explain the
reason.


\noindent The group $\; C_{N} \;$  acts from the right on the
manifold $\; \widetilde{01} \;$ with an action $\;T_{C_{N}}\;$
\begin{eqnarray} T_{C_{N}} & : &
\widetilde{01}\longrightarrow \widetilde{01} \nonumber \\
\label{jhhhdndnnmsk,ama,} ((g,h){\rm S}^{1}({\rm q}_{{\rm o}},{\rm
p}_{{\rm o}}))c & := & ((g,h)c){\rm S}^{1}({\rm q}_{{\rm o}},{\rm
p}_{{\rm o}})
 \end{eqnarray}
which is fixed point free. Since $\;C_{N}\;$ is finite and since
the action (\ref{jhhhdndnnmsk,ama,}) is fixed point free, the
coset space $\; 01 \equiv (({\rm S}{\rm L}(2,R) \times {\rm S}{\rm
L}(2,R)/{\rm S}^{1}({\rm q}_{{\rm o}},{\rm p}_{{\rm o}})) / C_{N}
\;$, which is the space of orbits of this action, inherits the
invariant measure of the space $\;  \widetilde{01} \equiv ({\rm
S}{\rm L}(2,R) \times {\rm S}{\rm L}(2,R))/{\rm S}^{1}({\rm
q}_{{\rm o}},{\rm p}_{{\rm o}}).\;$

\noindent This completes the necessary information
in order to construct  representations of $\; \mathcal H \mathcal
B
\;$ induced  from infinite little groups.

\section{Discussion}

\label{s5}

\noindent Two remarks are in order regarding the
representations of $\; \mathcal H \mathcal B
\;$ 
obtained by
the above construction
\begin{enumerate}
\item {As it was explained in \cite{Mel} \normalfont
the subgroup $\; 
L^{2}(\mathcal P,\lambda ,R)$ of 
$\; \mathcal H \mathcal B
= 
L^{2}(\mathcal P,\lambda ,R)
\bigcirc\!\!\;\!\!\!\;\;\;\;\;\;\;\;s \ _{T} \mathcal G \;$ is
topologised as a (pre) Hilbert space by using a natural measure
on $\; \mathcal P=P_{1}(R) \times P_{1}(R) \;$ and by introducing
a scalar product into $\;L^{2}(\mathcal P,\lambda ,R). \;$
If $\; R^{8} \;$ is endowed with the natural metric topology then
the group $\; \mathcal G=SL(2,R) \times SL(2,R), \;$ considered as
a subset of $\; R^{8}, \;$ inherits the induced topology on $\;\mathcal G.\;$
In the product topology of $\; 
L^{2}(\mathcal P,\lambda ,R) 
\times  \mathcal G
\;$ $\; \mathcal H \mathcal B
\;$ is a non$-$locally compact
group (the proof follows without substantial change Cantoni's
proof \cite{Cantoni}, see also \cite{mac3}). (In fact the subgroup
$\; 
L^{2}(\mathcal P,\lambda ,R),
\;$  and therefore the group $\; \mathcal H
\mathcal B
\;$ can be employed with many different topologies.
The Hilbert type topology employed here appears to describe
quantum mechanical systems in asymptotically flat space$-$times
\cite{Crampin2}).
Since in the Hilbert type topology
$\; \mathcal H \mathcal B
= 
L^{2}(\mathcal P,\lambda ,R)
\bigcirc\!\!\;\!\!\!\;\!\!\!\!s \  _{T} \mathcal G    \;$ is not
locally compact the theorems  dealing with the irreducibility of
the representations obtained by the above construction no longer
apply (see e.g. \cite{Mackey1}). However, 
it can be proved that the
induced representations obtained above \it {are} \normalfont irreducible. 
The proof
follows very closely the one given in \cite{mac5} for the case of
the original BMS group $\; B. \;$ }

\item { Here it is assumed that  $\; \mathcal H
\mathcal B
\;$ is equipped with the Hilbert topology. 
It is of outmost significance that it can be proved \cite{Mel} 
that in this topology $\; \mathcal H
\mathcal B
\;$ is a regular semi$-$direct$-$product. 
The proof follows the corresponding proof \cite{Piard1,Piard2} for the group $B$.
Regularity amounts  to the fact that \cite{Mackey} 
${L^{2}}^{'}(\mathcal P,\lambda ,R)$ can have no
equivalent classes of quasi$-$invariant measures $\; \mu \;$ such that the action of $\;
\mathcal G \;$ is strictly ergodic with respect to $\; \mu. \;$
When such measures $\; \mu \;$ do exist it can be proved \cite{Mackey} 
that 
an irreducible representation
of the group, with the semi$-$direct$-$product structure at hand, 
may be associated with each that is not equivalent to any of the IRs 
constructed by the Wigner$-$Mackey's inducing method.
In a different topology it is not known if $\; \mathcal H
\mathcal B
\;$ is a regular or irregular semi$-$direct$-$product.
Irregularity of $\; \mathcal H
\mathcal B$ in a topology different from the Hilbert 
topology would imply that there are IRs of $\; \mathcal H
\mathcal B$ that are not not equivalent to any of the
IRs obtained above by the inducing construction.
Strictly ergodic
actions are notoriously hard to deal with even in the locally
compact case. Indeed, for locally compact non$-$regular semi$-$direct
products, there is no known example for which all inequivalent
irreducibles arising from strictly ergodic actions have been
found. For the other 41 groups defined in
\cite{mac1} regularity has only been proved for $B$ (\cite{Piard1,Piard2})
when $B$ is equipped with the Hilbert topology.
Similar remarks apply to all of them regarding IRs arising from 
strictly ergodic actions in a given topology.

}

\end{enumerate}

\newpage

\begin{appendix}
\section{Fundamental Regions}

\label{s6}
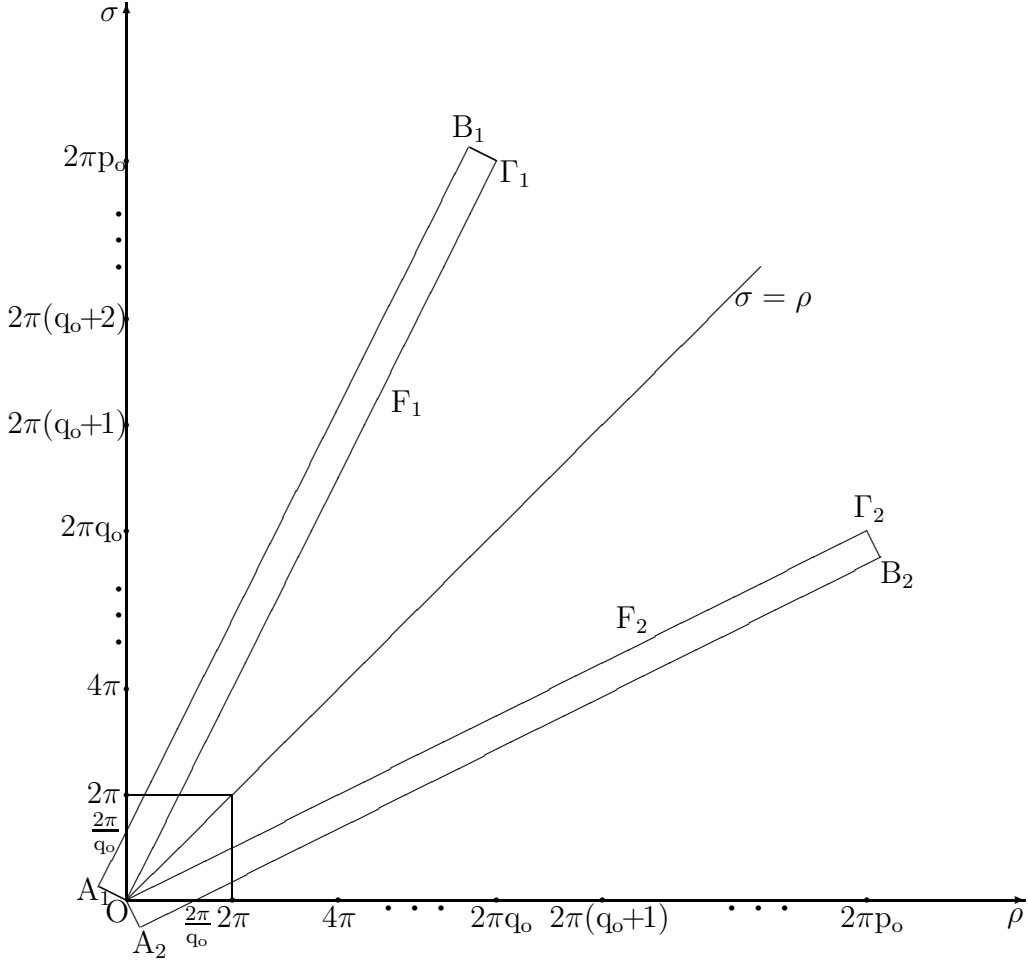
\begin{figure}[hp!]
\label{korasida} \setlength{\unitlength}{1pt}
\begin{center}
\begin{picture}(350,350)
\put(10,10){\vector(0,1){340}} \put(0,343){$\sigma$}
\put(10,10){\vector(1,0){340}} \put(343,2){$\rho$}
\put(10,50){\line(1,0){40}} \put(-4,33.5){$\frac{2\pi}{\rm q_{\rm
o}}$} \put(-5,47){$2\pi$}

\put(10,50){\circle*{2}} \put(-5,88){$4\pi$}
\put(10,90){\circle*{2}}


\put(110,195){${\rm F}_{1}$} \put(195,113){${\rm F}_{2}$}
\put(240,235){$ \sigma= \rho $}

\put(150,10){\circle*{2}} \put(190,10){\circle*{2}}

\put(140,0){$2 \pi {\rm q}_{\rm o}$} \put(170,0){$2 \pi ({\rm
q}_{\rm o}\!\!+\!\!1)$}

\put(290,10){\circle*{2}} \put(280,0){$2\pi{\rm p_{o}}$}

\put(50,10){\line(0,1){40}} \put(50,10){\circle*{2}}
\put(90,10){\circle*{2}} \put(84,0){$4 \pi $}

\put(109,7){\circle*{2}} \put(119,7){\circle*{2}}
\put(129,7){\circle*{2}}

\put(239,7){\circle*{2}} \put(249,7){\circle*{2}}
\put(259,7){\circle*{2}}

\put(31,-2){$\frac{2\pi}{\rm q_{o}}$}

\put(44,0){$2\pi$}

\put(10,10){\line(1,1){240}}

\put(10,10){\line(-2,1){10.3}} \put(10,10){\line(1,2){140}}
\put(-0.7,15){\line(1,2){140}} \put(150,290){\line(-2,1){10.3}}
\put(10,290){\circle*{2}} \put(-14,288){$2\pi{\rm p}_{\rm o}$}

\put(151,282){$\Gamma_{1}$} \put(133,299){${\rm B}_{1}$}
\put(-9,10){${\rm A }_{1}$} \put(2,2){O}

\put(285,155){$\Gamma_{2}$}
\put(10,150){\circle*{2}} \put(-15,148){$2\pi{\rm q}_{\rm o}$}
\put(10,190){\circle*{2}} \put(-35,187){$2\pi({\rm q}_{\rm
o}\!\!+\!\!1)$} \put(-35,227){$2\pi({\rm q}_{\rm o}\!\!+\!\!2)$}
\put(10,230){\circle*{2}} \put(7,108){\circle*{2}}
\put(7,118){\circle*{2}} \put(7,128){\circle*{2}}

\put(7,250){\circle*{2}} \put(7,260){\circle*{2}}
\put(7,270){\circle*{2}}

\put(295,130){${\rm B}_{2}$}

\put(12,-10){${\rm A}_{2}$}

\put(10,10){\line(1,-2){5}} \put(10,10){\line(2,1){280}}
\put(290,150){\line(1,-2){5}} \put(15.5,-0.1){\line(2,1){280}}

\end{picture}
\caption{\label{melas} \textbf{Fundamental regions  for the
families of groups
 $\mathbf{S^{1}_{(  q_ o,  p_ o)}}$,
 $\mathbf{S^{1}_{(  p_ o,  q_ o)}}$, \
where, $ \mathbf{p_ o> q_ o>0 \; \; {\bf or} \; \;
 p_ o= q_ o=1}$}.
The fundamental parallelograms ${\rm F}_{1}$ and ${\rm F}_{2}$
used in the text are depicted here. The co-ordinates of the
vertices of  ${\rm F}_{1}$ and ${\rm F}_{2}$ are: ${\rm A_{1}}:
\left(-\frac{{\rm q}}{1+{\rm q}^{2}} \frac {2 \pi}{{\rm q}_{\rm
o}} ,\frac {2 \pi}{1+{\rm q}^{2}} \frac {1}{{\rm q}_{\rm o}}
\right), \; {\rm B_{1}}: \left ( \frac {[{\rm q}_{\rm
o}^{2}(1+{\rm q}^{2})-{\rm q}]2 \pi} {(1+{\rm q}^{2}){\rm q}_{\rm
o}}, \frac {[{\rm q}_{\rm o}^{2}{\rm q}(1+{\rm q}^{2}) + 1]2 \pi}
{(1+{\rm q}^{2}){\rm q}_{\rm o}} \right ), \; \Gamma_{1} :(2 \pi
{\rm q}_{\rm o}, 2 \pi {\rm p}_{\rm o}), \; {\rm A_{2}} : \left (
\frac {2 \pi}{1+{\rm q}^{2}} \frac {1}{{\rm q}_{\rm o}},
-\frac{{\rm q}}{1+{\rm q}^{2}} \frac {2 \pi}{{\rm q}_{\rm o}}
\right), \; {\rm B_{2}} : \left ( \frac {[{\rm q}_{\rm o}^{2}{\rm
q}(1+{\rm q}^{2}) + 1]2 \pi} {(1+{\rm q}^{2}){\rm q}_{\rm o}} ,
\frac {[{\rm q}_{\rm o}^{2}(1+{\rm q}^{2})-{\rm q}]2 \pi} {(1+{\rm
q}^{2}){\rm q}_{\rm o}} , \right ), \; \Gamma_{2}: ( 2 \pi {\rm
p}_{\rm o}, 2 \pi {\rm q}_{\rm o}).\; {\rm The \; sides \;of \;
{\rm F}_{\rm 1} \;and \;{\rm  F}_{\rm 2} \; are \; given \; by :}
\; \; \quad {\rm O} \Gamma_{1} : \sigma= {\rm q} \rho , \; {\rm
A_{1}B_{1}}: \sigma = {\rm q} \rho + \frac { 2 \pi}{{\rm q}_{\rm
o}}, \; {\rm OA_{1}}: \sigma = - \frac {1}{{\rm q}} \rho , \; {\rm
B_{1}} \Gamma_{1} : \sigma = - \frac {1}{{\rm q}} \rho + \frac
{1+{\rm q}^{2}}{{\rm q}} 2 \pi {\rm q}_{\rm o}, \; {\rm
O}\Gamma_{2}: \rho={\rm q}\sigma, \; {\rm A_{2}B_{2}}: \rho ={\rm
q }\sigma + \frac {2 \pi }{ {\rm q}_{\rm o}}, \; {\rm OA_{2}}:
\rho = - \frac {1}{{\rm q}} \sigma , \; {\rm B_{2} \Gamma_{2}}:
\rho = - \frac {1}{\rm q} \sigma + \frac {1+{\rm q}^{2}}{\rm q} 2
\pi {\rm q}_{\rm o}.
 $}
\end{center}
\end{figure}


\begin{figure}[hp!]
\label{korasida} \setlength{\unitlength}{1pt}
\begin{center}
\begin{picture}(350,350)
\put(300,10){\vector(0,1){340}} \put(290,343){$\sigma$}
\put(10,10){\vector(1,0){340}} \put(-48,10){\line(1,0){60}}
\put(343,2){$\rho$} \put(260,50){\line(1,0){40}}
\put(301,33.5){$\frac{2\pi}{\rm q_{\rm o}}$} \put(301,47){$2\pi$}
\put(260,10){\line(0,1){40}} \put(246,0){$-2\pi$}
\put(260,10){\circle*{2}} \put(300,10){\line(-1,1){240}}
\put(80,235){$\sigma=-\rho$} \put(300,0){${\rm O}$}

\put(220,10){\circle*{2}} \put(206,0){$-4\pi$}
\put(143,0){$-2\pi{\rm q}_{\rm o}$} \put(160,10){\circle*{2}}

\put(120,10){\circle*{2}} \put(80,10){\circle*{2}}
\put(20,10){\circle*{2}}

\put(49,0){$-2\pi({\rm q}_{\rm o}\!\!+\!\!2)$}

\put(300,50){\circle*{2}} \put(304,88){$4\pi$}
\put(300,90){\circle*{2}}

\put(300,150){\circle*{2}} \put(304,148){$2\pi{\rm q}_{\rm o}$}
\put(300,190){\circle*{2}} \put(304,187){$2\pi({\rm q}_{\rm
o}\!\!+\!\!1)$} \put(304,227){$2\pi({\rm q}_{\rm o}\!\!+\!\!2)$}
\put(300,230){\circle*{2}} \put(304,108){\circle*{2}}
\put(304,118){\circle*{2}} \put(304,128){\circle*{2}}

\put(304,250){\circle*{2}} \put(304,260){\circle*{2}}
\put(304,270){\circle*{2}} \put(304,288){$2\pi{\rm p}_{\rm o}$}
\put(300,290){\circle*{2}}

\put(202,6){\circle*{2}} \put( 192,6){\circle*{2}}
\put(182,6){\circle*{2}}

\put(0,0){$-2\pi{\rm p}_{\rm o}$}

\put(30,6){\circle*{2}} \put( 40,6){\circle*{2}}
\put(50,6){\circle*{2}}


\put(300,10){\line(-1,2){139}} \put(300,10){\line(2,1){10}}
\put(160,289){\line(2,1){10}} \put(310,15){\line(-1,2){139}}

\put(311,12){${\rm A}_{3}$}

\put(173,291){${\rm B}_{3}$}

\put(187,200){${\rm F}_{3}$}

\put(148,286){${\rm \Gamma}_{3}$}

\put(300,10){\line(-2,1){281}} \put(294,0){\line(1,2){5}}
\put(294,0){\line(-2,1){280}} \put(14,141){\line(1,2){5}}

\put(105,112){${\rm F}_{4}$} \put(5,133){${\rm B}_{4}$}
\put(15,153){$\Gamma_{4}$} \put(287,-8){${\rm A}_{4}$}







\end{picture}
\caption{\label{bolos} \textbf{Fundamental regions  for the
families of groups
 $\mathbf{S^{1}_{(-  q_ o,  p_ o)}}$,
 $\mathbf{S^{1}_{(-  p_ o,  q_ o)}}$, \
where, $ \mathbf{p_ o> q_ o>0 \; \; {\bf or} \; \;
 p_ o= q_ o=1}$}.
The fundamental parallelograms ${\rm F}_{3}$ and ${\rm F}_{4}$
used in the text are depicted here. The co-ordinates of the
vertices of  ${\rm F}_{3}$ and ${\rm F}_{4}$ are: ${\rm A_{3}}:
\left(\frac{{\rm q}}{1+{\rm q}^{2}} \frac {2 \pi}{{\rm q}_{\rm o}}
,\frac {2 \pi}{1+{\rm q}^{2}} \frac {1}{{\rm q}_{\rm o}} \right),
\; {\rm B_{3}}: \left (- \frac {[{\rm q}_{\rm o}^{2}(1+{\rm
q}^{2})-{\rm q}]2 \pi} {(1+{\rm q}^{2}){\rm q}_{\rm o}}, \frac
{[{\rm q}_{\rm o}^{2}{\rm q}(1+{\rm q}^{2}) + 1]2 \pi} {(1+{\rm
q}^{2}){\rm q}_{\rm o}} \right ), \; \Gamma_{3} :(-2 \pi {\rm
q}_{\rm o}, 2 \pi {\rm p}_{\rm o}), \; {\rm A_{4}} : \left (-
\frac {2 \pi}{1+{\rm q}^{2}} \frac {1}{{\rm q}_{\rm o}},
-\frac{{\rm q}}{1+{\rm q}^{2}} \frac {2 \pi}{{\rm q}_{\rm o}}
\right), \; {\rm B_{4}} : \left (- \frac {[{\rm q}_{\rm o}^{2}{\rm
q}(1+{\rm q}^{2}) + 1]2 \pi} {(1+{\rm q}^{2}){\rm q}_{\rm o}} ,
\frac {[{\rm q}_{\rm o}^{2}(1+{\rm q}^{2})-{\rm q}]2 \pi} {(1+{\rm
q}^{2}){\rm q}_{\rm o}} , \right ), \; \Gamma_{4}: (- 2 \pi {\rm
p}_{\rm o}, 2 \pi {\rm q}_{\rm o}).\; {\rm The \; sides \;of \;
{\rm F}_{\rm 3} \;and \;{\rm  F}_{\rm 4} \; are \; given \; by :}
\; \; \quad {\rm O} \Gamma_{3} : \sigma= -{\rm q} \rho , \; {\rm
A_{3}B_{3}}: \sigma = -{\rm q} \rho + \frac { 2 \pi}{{\rm q}_{\rm
o}}, \; {\rm OA_{3}}: \sigma = - \frac {1}{{\rm q}} \rho , \; {\rm
B_{3}} \Gamma_{3} : \sigma =  \frac {1}{{\rm q}} \rho + \frac
{1+{\rm q}^{2}}{{\rm q}} 2 \pi {\rm q}_{\rm o}, \; {\rm
O}\Gamma_{4}: \rho=-{\rm q}\sigma, \; {\rm A_{4}B_{4}}: \rho
=-{\rm q }\sigma - \frac {2 \pi }{ {\rm p}_{\rm o}}, \; {\rm
OA_{4}}: \rho =  \frac {1}{{\rm q}} \sigma , \; {\rm B_{4}
\Gamma_{4}}: \rho =  \frac {1}{\rm q} \sigma - \frac {1+{\rm
q}^{2}}{\rm q} 2 \pi {\rm q}_{\rm o}.
 $}
\end{center}
\end{figure}
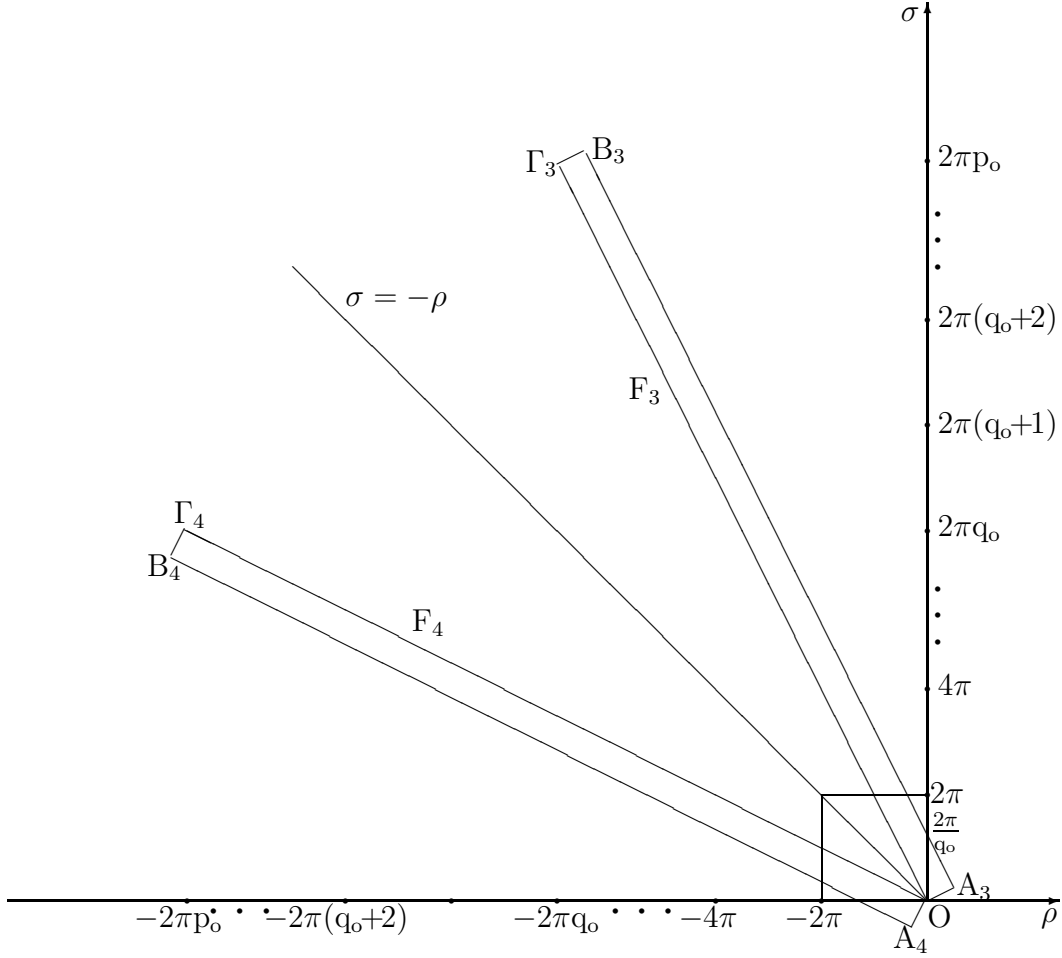


\section{$\mathcal G$$-$invariant measures}
\label{s7}

By following the corresponding analysis in \cite{mac3} we can
construct $\mathcal G$$-$invariant (and therefore 
$\;\mathcal G-{\rm quasi}-{\rm invariant})$ measures on the orbits $ \mathcal G
/L_{\phi_{\rm{o}}} \approx \mathcal G \phi_{\rm{o}} .$ Firstly, we need the
following Theorem (\cite{Helgason}, p.369):

\noindent {\bf Theorem} {\em Let $\mathcal G$ be a Lie group, $H$
a closed subgroup. The relation \be |\det {\rm A}{\rm d}_{\mathcal
G}(h)|=|\det {\rm A}{\rm d}_{H}(h)| \quad (h\in H) \ee is a
necessary and sufficient condition for the existence of a positive
$\mathcal G$$-$invariant measure ${\rm d}\mu_{H}$ on $\mathcal G/H$,
which is unique up to a constant factor, and satisfies \be
\label{skiros} \int_{\mathcal G}f(g){\rm d}\mu(g)=\int_{\mathcal
G/H} \left( \int_{H}f(gh){\rm d}\mu(h) \right) {\rm d}\mu_{H}, \ee
where ${\rm d}\mu(g)$ \  and \ ${\rm d}\mu(h)$ are, respectively,
suitably normalized invariant measures on $\mathcal G$ and $H$.}

\noindent Here ${\rm A}{\rm d}_{\mathcal G}$ denotes the adjoint
representation of the group $\mathcal G$ and $f$ is any continuous
function of compact support on $\mathcal G$. Since any function on
$\mathcal G$ constant on $H$ cosets may be regarded as a function
on $\mathcal G/H$, this defines ${\rm d}\mu_{H}$ on $\mathcal
G/H$. (The space of continuous, compact support functions may be
replaced by the space of integrable functions needed  in the text
by the usual completion process). Thus, in each case, it is
sufficient to verify the condition on the moduli of the
determinants, and to provide an $H$ invariant measure ${\rm
d}\mu(h)$ on $H$.

\indent Since every $H=L_{\phi_{\rm{o}}}$ is abelian, ${\rm A}{\rm
d}_{L_{\phi_{\rm{o}}}}(\varepsilon)$ is the identity operator, so that
$|\det {\rm A}{\rm d}_{L_{\phi_{\rm{o}}}}(\varepsilon)|=1$ \ for all \
$L_{\phi_{\rm{o}}}$ \ and for all \ $\varepsilon \in L_{\phi_{\rm{o}}}$. It must now
be shown that $|det Ad_{\mathcal G}(\varepsilon)|=1$.\  A basis
for the Lie Algebra of $\mathcal G$ is given by the generators \be
\Omega_{i}=\left [ \begin{array}{ll}
A_{i} & {\bf 0} \\
{\bf 0} & {\bf 0}
\end{array}
\right ], \ \Omega_{3+i}=\left [ \begin{array}{ll}
{\bf 0} & {\bf 0} \\
{\bf 0} & A_{i}
\end{array}
\right ], \ee where ${\bf 0}$ is the $2\times 2$ zero-matrix and
$A_{i},i=1,2,3$ \be A_{1}=\left [ \begin{array}{lr}
0 & 1 \\
1 & 0
\end{array}
\right ], \ A_{2}=\left [ \begin{array}{lr}
1& 0 \\
0 & -1
\end{array}
\right ], \ A_{3}=\left [ \begin{array}{rr}
0& 1 \\
-1 & 0
\end{array}
\right ], \ee is a basis for the Lie Algebra of ${\rm S}{\rm
L}(2,{\rm R})$.\ Taking \be \varepsilon =  \left[
\begin{array}{cccc}
\cos({\rm p}_{\rm o}\theta) & \sin({\rm p}_{\rm o}\theta) & 0 & 0  \\
                      \!\!\!\! -\sin({\rm p}_{\rm o}\theta) &
\cos({\rm p}_{\rm o}\theta) & 0 & 0  \\
                        0 & 0 & \cos({\rm q}_{\rm o}\theta) &
\sin({\rm q}_{\rm o}\theta)  \\
                        0 & 0 & \!\!\!\!\!-\sin({\rm q}_{\rm o}\theta) &
\cos({\rm q}_{\rm o}\theta)
\end{array}
\right], \ee a straightforward calculation shows \be {\rm A}{\rm
d}_{\mathcal G}(\varepsilon)=  \left[
\begin{array}{cccccc}
\cos(2{\rm p}_{\rm o}\theta) & \sin(2{\rm p}_{\rm o}\theta) & 0 & 0 & 0 & 0 \\
\!\!\!\!-\sin(2{\rm p}_{\rm o}\theta) & \cos(2{\rm p}_{\rm o}\theta) & 0 & 0 & 0 & 0 \\
0 & 0 &  1 & 0 & 0 & 0 \\
0 & 0 &  0 & \cos(2{\rm q}_{\rm o}\theta) & \sin(2{\rm q}_{\rm o}\theta) & 0 \\
0 & 0 &  0 & \!\!\!\!-\sin(2{\rm q}_{\rm o}\theta) & \cos(2{\rm q}_{\rm o}\theta) & 0 \\
0 & 0 &  0 & 0 & 0 & 1
\end{array}
\right ] \ee Evidently $\det {\rm A}{\rm d}_{\mathcal
G}(\varepsilon)=1$\ for all $\varepsilon \in L_{\phi_{\rm{o}}}$ \ and all
$L_{\phi_{\rm{o}}}$,\ so that the required condition is satisfied. When
${\rm p}_{{\rm o}},{\rm q}_{{\rm o}}$ are relatively prime, the
limits of integration on  $\int_{H}f(gh){\rm d}\mu(h)$
(Eq.(\ref{skiros})) are 0 and $2\pi$. A normalised $\; L_{\phi_{\rm{o}}}
\;$$-$invariant measure on the one$-$dimensional little groups
$L_{\phi_{\rm{o}}}$ is $\frac{1}{2 \pi}{\rm d}\theta$. Hence the $\mathcal
G$$-$invariant measures  on $\mathcal G/L_{\phi_{\rm{o}}}$ are  given.

Next, the homeomorphism types of the orbits will be given. Recall
that every $g \in {\rm S}{\rm L}(2,{\rm R})$ has a unique `polar'
decomposition \be g= \kappa u, \ee where $\kappa$ is a positive
semi$-$definite symmetric matrix, and $u \in {\rm S}{\rm O}(2)$. In
the usual metric topology of ${\rm S}{\rm L}(2,{\rm R})$, it is
easily shown that the set of $\kappa$'s forms a topological space
homeomorphic to $R^{2}$. Since $\;{\rm S}{\rm O}(2)\;$ and $\;{\rm
S}^{1}_{( {\rm q}_{\rm o}, {\rm p}_{\rm o})} \;$ have the
topological structure of a circle, we conclude that the orbit of
the ${\rm S}^{1}_{( {\rm q}_{\rm o}, {\rm p}_{\rm o})}$$-$dual
action through a point $\phi$ in $A'(\N)$ is homeomorphic to: \be
R^{4} \times {\rm S}^{1}
. \ee
\end{appendix}

\bibliographystyle{plain}
{}

\end{document}